\documentclass[a4paper]{amsart}

\usepackage{amscd,amsthm,amsfonts,amssymb,amsmath}
\usepackage[colorlinks=true]{hyperref}
\usepackage{fullpage}
\usepackage[all]{xy}

    \newcommand{\BA}{{\mathbb {A}}} \newcommand{\BB}{{\mathbb {B}}}
    \newcommand{\BC}{{\mathbb {C}}} 
     \newcommand{\BF}{{\mathbb {F}}}

     \newcommand{\BP}{{\mathbb {P}}}
    \newcommand{\BQ}{{\mathbb {Q}}} \newcommand{\BR}{{\mathbb {R}}}
    \newcommand{\BS}{{\mathbb {S}}} \newcommand{\BT}{{\mathbb {T}}}

     \newcommand{\BZ}{{\mathbb {Z}}}

     \newcommand{\CF}{{\mathcal {F}}}
     \newcommand{\CH}{{\mathcal {H}}}

    \newcommand{\CO}{{\mathcal {O}}} \newcommand{\CP}{{\mathcal {P}}}
     \newcommand{\CR}{{\mathcal {R}}}

     \newcommand{\RN}{{\mathrm {N}}}

     \newcommand{\fp}{{\mathfrak{p}}}

    \newcommand{\Aut}{{\mathrm{Aut}}}

    \newcommand{\End}{{\mathrm{End}}}

    \newcommand{\Gal}{{\mathrm{Gal}}} \newcommand{\GL}{{\mathrm{GL}}}
    
    \newcommand{\Hom}{{\mathrm{Hom}}}

    \newcommand{\Isom}{{\mathrm{Isom}}}

    \newcommand{\ord}{{\mathrm{ord}}} \newcommand{\rank}{{\mathrm{rank}}}
     \newcommand{\Pic}{\mathrm{Pic}}

    \renewcommand{\mod}{\ \mathrm{mod}\ }

    \newcommand{\Res}{{\mathrm{Res}}}
    \newcommand{\Sel}{{\mathrm{Sel}}}

    \newcommand{\SL}{{\mathrm{SL}}}
     
    \newcommand{\SO}{{\mathrm{SO}}}

    \newcommand{\tr}{{\mathrm{tr}}}\newcommand{\tor}{{\mathrm{tor}}}
    \newcommand{\RTr}{{\mathrm{Tr}}}
    
    \newcommand{\Vol}{{\mathrm{Vol}}}

        \newcommand{\Tr}{\mathrm{Tr}}

\newcommand{\matrixx}[4]{\begin{pmatrix}
#1 & #2 \\ #3 & #4
\end{pmatrix} }        % 2*2 matrix

    \font\cyr=wncyr10

    \newcommand{\Sha}{\hbox{\cyr X}}
    \newcommand{\wh}{\widehat}

    \newcommand{\ov}{\overline}

    \newcommand{\sk}{\medskip}
    \newcommand{\lra}{\longrightarrow}
    \newcommand{\ra}{\rightarrow} 
    \newcommand{\lto}{\longmapsto}\newcommand{\bs}{\backslash}
    
    \newcommand{\s}{\sk\noindent}

    \theoremstyle{plain}
    \newtheorem{thm}{Theorem}[section] \newtheorem{cor}[thm]{Corollary}
    \newtheorem{lem}[thm]{Lemma}  \newtheorem{prop}[thm]{Proposition}
     \newtheorem{defn}[thm]{Definition}

\theoremstyle{remark} 
\theoremstyle{remark} 
\theoremstyle{remark} 

    \newcommand{\Neron}{N\'{e}ron~}

    \numberwithin{equation}{section}

\begin{document}

\title{Cube Sum Problem and an Explicit Gross-Zagier Formula}

\author{Li Cai, Jie Shu,  and Ye Tian}

\address{Li Cai: Mathematical Sciences Center, Tsinghua University, Beijing
100084} \email{lcai@math.tsinghua.edu.cn}

\address{Jie Shu: Academy of Mathematics and Systems
Science, Morningside center of Mathematics, Chinese Academy of
Sciences, Beijing 100190} \email{shujie09@mails.gucas.ac.cn}

\address{Ye Tian: Academy of Mathematics and Systems
Science, Morningside center of Mathematics, Chinese Academy of
Sciences, Beijing 100190} \email{ytian@math.ac.cn}

\thanks{Li Cai was supported by the Special Financial Grant from the China Postdoctoral Science Foundation 2014T70067; 
Ye Tian was supported by  the 973 Program 2013CB834202 and NSFC grants 11325106, 11321101, and 11031004}

%\keywords{Elliptic Curve}

\maketitle

\tableofcontents

\section{Introduction and Main Results}

We call a nonzero rational number a {\em cube sum} if it is of form $a^3+b^3$ with $a,b\in \BQ^\times$. Similar to Heegner's work on congruent number problem, Satg\'e \cite{Satge} showed that if  $p\equiv 2\mod 9$ (resp. $p\equiv 5\mod 9$) is a prime, then $2p$ (resp. $2p^2$) is a cube sum. Lieman \cite{Lieman} showed there are infinitely many cube-free integers which are not cube sums. Sylvester \cite{Syl} showed there are infinitely many integers with at most $2$ distinct prime factors which are not cube sums. The following is one of main results in this paper.
\begin{thm}\label{th:1}
For any odd integer $k\geq 1$, there exist infinitely many cube-free odd integers $n$ with exactly $k$ distinct prime factors such that $2n$ is a cube sum (resp. not a cube sum).
\end{thm}
\par For any $n\in \BQ^\times$, let $C^{(n)}$ be the elliptic curve over $\BQ$ defined by the equation $x^3+y^3=2n$.  Note that the torsion part of the Mordell-Weil group $C^{(n)}(\BQ)$ is trivial unless $2n$ is a cube, which is not a cube sum by our convention.  Then $2n$ is a cube sum if and only if the rank of the Mordell-Weil group $C^{(n)}(\BQ)$ is positive.

For an odd prime $p\equiv 2,5\mod 9$, denote $p^*=p^{\pm 1}\equiv 2\mod 9$. To prove that $p^*$ is a cube sum,  Satg\'e \cite{Satge} constructed a non-trivial Heegner point on $C^{(p^*)}$ (see also Dasgupta and Voight \cite{DV}.)
In this paper, we give a similar construction of Heegner point on $C^{(p^*)}$ and
relate its height to some special $L$-value.
Together with work of Kolyvagin \cite{K}, Perrin-Riou \cite{Perrin-Riou} and Kobayashi \cite{Kobayashi},
we have the following result on the Birch and Swinnerton-Dyer conjecture for $C^{(p^*)}$ and $C^{(p^{*-1})}$.
\begin{thm}\label{th:2}
Let $p\equiv 2,5\mod 9$ be an odd prime number. Then $2p^*$ is a cube sum. Moreover,
\begin{enumerate}
\item
$\ord_{s=1} L(s,C^{(p^*)})=\rank_\BZ\,C^{(p^*)}(\BQ)=1$  and $\ord_{s=1} L(s,C^{(p^{*-1})})=\rank_\BZ\,C^{(p^{*-1})}(\BQ)=0$.
\item  The Tate-Shafarevich groups $\Sha(C^{(p)})$ and $\Sha(C^{(p^{-1})})$ are finite, and, for any prime $\ell\nmid 2p$,
	the  $\ell$-part of   $\#\Sha(C^{(p)})\cdot \#\Sha(C^{(p^{-1})})$
	is as predicted by the Birch and Swinnerton-Dyer conjecture for $C^{(p)}$ and $C^{(p^{-1})}$.
\end{enumerate}
\end{thm}

\par

In this paper,  we give also a general construction of Heegner point and obtain an  explicit Gross-Zagier formula,   which is
a variant of our previous work \cite{CST} and is used to prove Theorem $\ref{th:2}$.

\par Let $\phi$ be a newform of weight $2$, level $\Gamma_0(N)$, with
Fourier expansion $\phi = \sum_{n=1}^\infty a_n q^n$ normalized
such that $a_1 = 1$. Let $K$ be an imaginary quadratic field
of discriminant $D$ with $\CO_K$ its
ring of integers. For an positive integer $c$, let
$\CO_c = \BZ + c\CO_K$ be the order of $K$ of conductor $c$. Denote by
$H_c$ the ring class field of $K$ which is the abelian extension over $K$
characterized by the property that the Artin map induces an isomorphism
$\Pic(\CO_c) \stackrel{\sim}{\ra} \Gal(H_c/K)$. Here $\Pic(\CO_c)$ is
the Picard group for invertible (fractional) $\CO_c$-ideals.
Let $\chi:\Gal(H_c/K) \ra \BC^\times$ be a primitive character of
ring class field.  Let $L(s,\phi,\chi)$ be the
Rankin-Selberg $L$-function of $\phi$ and $\chi$ which is defined by
an Euler product over primes $p$
\[L(s,\phi,\chi) = \prod_{p<\infty} L_p(s,\phi,\chi).\]
We refer to \cite{GJ} for the general definition of the local factors $L_p(s,\phi,\chi)$.
Denote by $S$ the set of finite primes $p|(N,Dc)$ such that if
$p || N$ then $p|c$. For any $p \not \in S$, the local
factor $L_p(s,\phi,\chi) = \prod_{1 \leq i,j \leq 2} (1-\alpha_{p,i}\beta_{p,j}p^{-s})^{-1}$
where $\alpha_{p,i}$ and $\beta_{p,j}$ are local parameters of $L$-series $L(s,\phi)$ and $L(s,\chi)$:
\[\prod_{i = 1,2} (1-\alpha_{p,i}p^{-s}) = L_p(s,\phi)^{-1} =
	\begin{cases}
		1- a_pp^{-s} + p^{1-2s}, \quad &\text{if } p \nmid N; \\
		1- a_pp^{-s}, \quad &\text{if } p || N;\\
		1, \quad &\text{if } p^2|N,
\end{cases}\]
and
\[
	\prod_{j = 1,2} (1-\beta_{p,j}p^{-s}) = L_p(s,\chi)^{-1} =
	\begin{cases}
		\prod_{\fp|p} (1-\chi([\fp])N\fp^{-s}), \quad &\text{if } p\nmid c; \\
		1, \quad &\text{if } p|c,
	\end{cases}
	\]
where $[\fp] \in \Pic(\CO_c)$ is the ideal class of the  prime $\fp|p$.
The $L$-function $L(s,\phi,\chi)$ admits an analytically continuation to the whole $s$-plane
with the functional equation
\[G_2(s)^2L(s,\phi,\chi) = \epsilon(s,\phi,\chi)G_2(2-s)^2L(2-s,\phi,\chi)\]
where $G_2(s) := 2(2\pi)^{-s}\Gamma(s)$ and the root number
$\epsilon(\phi,\chi) := \epsilon(1,\phi,\chi)  = \pm 1$.

Firstly, assume that:
\begin{enumerate}
	\item $(c,N) = 1$, and if $p|(N,D)$ then $\ord_p(N) = 1$ (in particular,
	      $a_p = \pm 1$).
	\item The set  $\Sigma$  has even cardinality, where $\Sigma$ consists of prime factors $p$ of $N$ satisfying
		\begin{itemize}
			\item either $p$ is inert in $K$ and $\ord_p(N)$ is odd,
			\item or  $p$ is ramified in $K$ and $\chi([\fp]) = a_p$
				where $[\fp]\in \Pic(\CO_c)$ is the ideal class of the prime ideal $\fp$ above $p$.
		\end{itemize}
	
\end{enumerate}
The above assumption implies (see Lemma 3.1 in \cite{CST}) that the root number $\epsilon(\phi,\chi)$ of $L(s,\phi,\chi)$ is $-1$.
Let $B$ be the indefinite quaternion algebra defined over $\BQ$
ramified exactly at the set $\Sigma$. There is
then an $\BQ$-embedding from $K$ to $B$ and we shall fix one from now on.
Let $R$ be an order
in $B$ with discriminant $N$ such that $R \cap K = \CO_c$. Such $R$ exists.
%It can be constructed locally as following.
%For a finite place $p$, if $M$ is a $\BZ$-module, denote by $M_p := M \otimes_\BZ \BZ_p$.
%For each $p$, there is a local order $\tilde{R}_p$ of $B_p$ with discriminant $N\BZ_p$
%such that $\tilde{R}_p \cap K_p = \CO_{c,p}$ and such order is unique up to
%$K_p^\times$-conjugacy (See also Lemma 3.3 in \cite{CST}). Moreover, we may assume that
%for almost all $p$, $\tilde{R}_p = R_{0,p}$ for some maximal order $R_0$ of $B$. Then
%$R := B \cap (\prod_p \tilde{R}_p)$ is an order of $B$ we want.
%Moreover, if denote by $N_B$ the reduced norm for $B$ (also for each $B_p$),
%then for any finite $p$, $N_B(R_p^\times) = \BZ_p^\times$. To see this, note that
%if $p \not| (D,N)$, then either $R_p$ is maximal or contains $\CO_{K,p}$ and the image of
%$R_p^\times$ under $N_B$ is $\BZ_p^\times$. On the other hand, if $p| (D,N)$, then
%$\ord_p(N) = 1$ by our assumption which implies that $R_p$ is either an Eichler order with
%discriminant $p\BZ_p$ if $B_p$ is split or the maximal order of $B_p$ if $B_p$ is nonsplit.
Let $\RN_B$ be the reduced norm of $B$. Put
\[\Gamma = \left\{ \gamma \in R^\times | \RN_B(\gamma) = 1 \right\}.\]
%\Gamma = B_+^\times \cap \wh{R}^\times
The group $\Gamma$ is viewed as a subgroup of $\SL_2(\BR)$ via the embedding
$B \hookrightarrow B\otimes_\BQ \BR \cong M_2(\BR)$ which
is a Fuchsian group of the first kind (See also \cite{Miyake}, Theorem 5.2.13).
Let $X_\Gamma$ be the Shimura curve over $\BQ$  associated to $B$ with level $\Gamma$.
The curve $X_\Gamma$ is a geometrically connected projective curve which
has the complex uniformization
\[X_\Gamma(\BC) = \Gamma \bs \CH \cup \{\mathrm{cusps}\}\]
where $\CH$ is the upper half plane with the linear fractional
action of $\Gamma$. The set of cusps is non-empty if and only if
$\Sigma = \emptyset$.
%Note that by our assumption, if we denote $U = \wh{R}^\times$,
%then $N_B(U) = \wh{\BZ}^\times$ and $\wh{B}^\times = B_+^\times U$.
Let $h_0$ be the unique point in $\CH$ fixed by $K^\times$ and
denote by $P$ the point on $X_\Gamma(\BC)$ represented by $h_0$.
Then by Shimura's reciprocity law, as $R \cap K = \CO_c$,
$P$ is defined over the ring class field $H_c$ over $K$ of conductor $c$.

Let $L\in \Pic(X_\Gamma)_\BQ$ be the Hodge class. It is the line bundle whose global sections are holomorphic modular forms of weight two, i.e.
\[L=\omega_{X_\Gamma/\BQ} + \sum_{x \in X_\Gamma(\bar{\BQ})} (1-e_x^{-1})x.\]
Here $\omega_{X_\Gamma/\BQ}$ is the canonical bundle of $X_\Gamma$,
$e_x$ is the ramification index of $x$ in the complex uniformization of
$X_\Gamma$, i.e. for a cusp $x$, $e_x = \infty$ and for a non-cusp $x$,
$e_x$ is half of the order of stabilizers of $x$ in $\Gamma$. Let $\xi=L/\deg L\in \Pic(X_\Gamma)_\BQ$ be the normalized Hodge class on $X_\Gamma$ (of degree one). 
 Here $\deg L$ denotes the degree of $L$, which is the volume of $X_\Gamma(\BC)$ with respect to the measure $dxdy/2\pi y^2$.  Let $J_\Gamma$ be the Jacobian of $X_\Gamma$. The degree one class $\xi$ defines a morphism 
 $$[\quad ]: \ X_\Gamma(\bar{\BQ})\lra J_\Gamma(\bar{\BQ})_\BQ, \qquad x \lto [x-\xi].$$  Consider the point
$$P_\chi:=\sum_{\sigma\in \Gal(H_c/K)} [P-\xi]^{\sigma} \otimes \chi(\sigma)\in J_\Gamma(H_c)_\BC := J_\Gamma(H_c)_\BQ\otimes_\BQ \BC.$$
Let $\BT\subset \End_\BQ(J_\Gamma)$ be the sub-algebra generated over $\BZ$ by Hecke correspondences $T_\ell$, for $\ell$ prime to $N$. Let $P_\chi^\phi$ be the projection of $P_\chi\in J_\Gamma(H_c)_\BC$ to the eigenspace with the 
eigen character $\BT\ra \BC$ defined by $T_\ell \mapsto a_\ell$.

\begin{thm}
Under assumptions (1) and (2) above, we have the following identity
  \[L'(1,\phi,\chi) = 2^{-\mu(N, D)}\cdot
    \frac{8\pi^2(\phi,\phi)_{\Gamma_0(N)}}{u^2\sqrt{|Dc^2|}}\cdot
  \wh{h}_K(P_\chi^\phi),\]
  where $\mu(N, D)$ is the number of common prime factors of $N$ and $D$, $u = [\CO_c^\times: \BZ^\times]$
  and
  \[(\phi,\phi)_{\Gamma_0(N)} = \iint_{\Gamma_0(N) \bs \CH}
  |\phi(x+iy)|^2dxdy,\]
  and $\wh{h}_K$ is the \Neron-Tate height on $J_\Gamma(H_c)_\BC$ over $K$. 
\end{thm}

Let $A$ be an abelian variety over $\BQ$ associated to $\phi$ in the sense that for all primes $p$,
$$L_p(s, A)=\prod_{\sigma: \BQ_\phi\hookrightarrow \BC} L(s, \phi^{\sigma}),$$
where $\BQ_\phi\subset \BC$ is the subfield generated over $\BQ$ by all Fourier coefficients of $\phi$. Then  $M:=\End^0(A)$ is a number field over $\BQ$ of degree $\dim A$. There exists a unique embedding $\iota: M\hookrightarrow \BC$  satisfying the following property:
for any $p \nmid N$,  if we take $\ell \not= p$ a prime inert in $M$, then under $\iota$, the trace of the
$M_\ell$-linear map defined by the action of the Frobelinus element  at $p$ on the
$\ell$-adic Tate module $V_\ell(A) := \lim_{n} A(\bar{\BQ})[\ell^n] \otimes_{\BZ_\ell} \BQ_\ell$ equals $a_p$.

 \begin{defn} A non-constant morphism $f: X_\Gamma \ra A$ over $\BQ$ is called 
	 an admissible modular parametrization if there is an integral multiple of $\xi$  
	 represented by a divisor $\sum n_i x_i$ with integral coefficients $n_i$ such that  
	 $\sum n_i f(x_i)$ is equal to the identity element in $A(\BQ)$.
\end{defn}
\begin{prop}[\cite{CST} Proposition 3.8]
There exists an admissible modular parametrization $f: X_\Gamma \ra A$. 
Moreover, if $f'$ is another admissible modular parametrization, then there exist nonzero integers
$n,n'$ such that $nf = n'f'$.
\end{prop}

%Let $p$ be a prime dividing $N$ and nonsplit in $K$. Let $\varpi_p$ be a uniformizer
%of $K_p$. Then $R_p$ is normalized by $\varpi_p$ (See also Lemma 3.4 in \cite{CST}). The Hecke action of $\varpi_p$
%is an automorphism on $X_\Gamma$ over $\BQ$. On the other hand, as $\chi$ is unramified at $p$, $\chi([\fp]) \in \{\pm 1\}$
%where $[\fp] \in \Pic(\CO_c)$ is the ideal class of the prime $\fp$ over $p$. In this case,
%the modular parametrization defined by the composition of $f$ with the Hecke action of $\varpi_p$ equals $\chi([\fp])f$.

By the above proposition, we may take an admissible modular parametrization $f_1: X_\Gamma \ra A$.
Let $A^\vee$ be the dual abelian variety of $A$. Since $A^\vee$ is isogeny to $A$, there is also an admissible modular
parametriation $f_2: X_\Gamma \ra A^\vee$. View $f_1 \in \Hom(J_\Gamma,A)$ and
$f_2 \in \Hom(J_\Gamma,A^\vee)$. Denote
\[(f_1,f_2)_\Gamma := f_1 \circ f_2^{\vee} \in M \stackrel{\iota}{\hookrightarrow} \BC\]
where $f_2^{\vee}: A \ra J_\Gamma$ is the dual of $f_2$ composed with
the canonical isomorphism $J_\Gamma^\vee \cong J_\Gamma$. If $A$ is
an elliptic curve and identify $A^\vee$ with $A$ canonically, then for
any morphism $f:X_\Gamma \ra A$, we have $(f,f)_\Gamma = \deg f$, the
degree of $f$.
Let
\[P_\chi(f_1) = \sum_{\sigma \in \Gal(H_c/K)} f_1(P)^\sigma \chi(\sigma) \in A(H_c)_\BC := A(H_c)_\BQ \otimes_{(M,\iota)} \BC.\]
Similarly, define $P_{\chi^{-1}}(f_2) \in A^\vee(H_c)_\BC$.

The usual theory of N\'eron-Tate height over $K$ gives a $\BQ$-bilinear non-degenerated pairing
(See also \cite{YZZ}, \S 7.1.1)
\[\langle \cdot,\cdot \rangle_K: A(\bar{K})_\BQ \times A^\vee(\bar{K})_\BQ \lra \BR.\]
The field $M$ acts on $A(\bar{K})_\BQ$ by definition, and acts on
$A^\vee(\bar{K})_\BQ$ through the duality.
By the adjoint property of the height pairing, the Neron-Tate height pairing $\langle \cdot,\cdot \rangle_K$
descends to  a $\BQ$-linear map
\[\langle \cdot,\cdot \rangle_K: A(\bar{K})_\BQ \otimes_{M} A^\vee(\bar{K})_\BQ \lra \BR.\]
Denote by $V$ the $M$-module $A(\bar{K})_\BQ \otimes_{M} A^\vee(\bar{K})_\BQ$ for short. The
$\BQ$-linear map $\langle \cdot,\cdot \rangle_K: V \ra \BR$ induces
a $\BC$-linear map $V \otimes_\BQ \BC \ra \BC$. Note that
\[V \otimes_\BQ \BC = V \otimes_M (M \otimes_\BQ \BC) =
\bigoplus_{\iota': M \hookrightarrow \BC} V \otimes_{(M,\iota')} \BC\]
where $\iota'$ runs over all the embeddings from $M$ to $\BC$.
The above $\BC$-linear map induces a $\BC$-linear map $V \otimes_{(M,\iota)} \BC \ra \BC$
which gives an $M$-linear map $V \ra \BC$
with $M$ acting on $\BC$ via $\iota$. This further induces the following $M$-linear map 
with the action of $M$ on $\BC$ given by $\iota$
\[\langle \cdot,\cdot \rangle_K: A(\bar{K})_\BC
\otimes_{M} A^\vee(\bar{K})_\BC \lra \BC\]
where  $A(\bar{K})_\BC := A(\bar{K})_\BQ \otimes_{(M,\iota)} \BC$ and
$A^\vee(\bar{K})_\BC := A^\vee(\bar{K})_\BQ \otimes_{(M,\iota)} \BC$.

%Let $\langle , \rangle_{K,M}:
%A(\bar{K})_\BQ \otimes_{M} A^\vee(\bar{K})_\BQ \ra M\otimes_\BQ \BR$ be the unique $M$-bilinear
%      pairing such that $\langle , \rangle_K = \tr_{M\otimes \BR/\BR} \langle , \rangle_{K,M}$. This
%      induces a $L$-linear Neron-Tate pairing over $K$:
%      \[\left\langle \cdot,\cdot \right\rangle_{K,L}:
%	      (A(\bar{K})_\BQ \otimes_{M} L) \otimes (A^\vee(\bar{K})_\BQ \otimes_{M} L)
%      \lra L\otimes_\BQ \BR.\]

\begin{thm}\label{gz-disjoint}
  Under assumptions (1) and (2) above, we have the following identity
  \[L'(1,\phi,\chi) = 2^{-\mu(N, D)}\cdot
    \frac{8\pi^2(\phi,\phi)_{\Gamma_0(N)}}{u^2\sqrt{|Dc^2|}}\cdot
  \frac{\langle P_\chi(f_1),P_{\chi^{-1}}(f_2)\rangle_{K}}{(f_1,f_2)_\Gamma}.\]
  Here, $\mu(N, D)$ is the number of common prime factors of $N$ and $D$, $u = [\CO_c^\times: \BZ^\times]$
  and
  \[(\phi,\phi)_{\Gamma_0(N)} = \iint_{\Gamma_0(N) \bs \CH}
  |\phi(x+iy)|^2dxdy.\]
\end{thm}

However, for cube sum problem, the assumption $(c,N) = 1$ does not hold and we need to
carefully choose the test vector $f$ at any place $p|(c,N)$. Based on the work of
Yuan-Zhang-Zhang \cite{YZZ}, this
is in fact a problem of harmonic analysis on local $p$-adic representations. 
Let $B$ be a quaternion algebra over $\BQ_p$ for some $p$ with a
quadratic sub $\BQ_p$-algebra, say $K$. Denote by $G = B^\times$. Let $\pi$ be an irreducible
smooth unitary representation on $G$ which is  of infinite dimension if
$B$ is split. Assume that the central character of $\pi$ is trivial. Let $\chi$ be a
character on $K^\times$ such that $\chi|_{\BQ_p^\times} = 1$. Consider the functional space 
\[\CP(\pi,\chi) := \Hom_{K^\times}(\pi,\chi^{-1}).\] 
In general, its dimension is less than one and here we assume the dimension is one. 
A nonzero vector $f \in \pi$ is called a {\em test vector} for the pair $(\pi,\chi)$ provided that
$\ell(f) \not=0$ for any nonzero functional $\ell \in \CP(\pi,\chi)$. For a
nonzero $f \in \pi$, consider the following toric integral
\[\beta^0(f) = \int_{\BQ_p^\times \bs K^\times} \frac{(\pi(t)f,f)}{(f,f)} \chi(t)dt\]
where $(\cdot,\cdot)$ is any nonzero invariant Hermitian pairing on $\pi$ and $dt$ is any
Haar measure on $\BQ_p^\times \bs K^\times$. This integral is absolutely convergent. 
Moreover, $f$ is a test vector for $(\pi,\chi)$ if and only if $\beta^0(f) \not=0$.

Denote by $n$ the conductor of $\pi$, that is, the minimal non-negative integer
such that the invariant subspace of $\pi$ under $U_0(n)$ is nonzero. Here
$U_0(n) = R_0(n)^\times$ with 
$R_0(n) = 
\begin{pmatrix}
	\BZ_p & \BZ_p \\
	p^n\BZ_p & \BZ_p
\end{pmatrix}$. Denote by $c$ the conductor of $\chi$, that is, the minimal
non-negative integer such that $\chi$ is trivial on the units of 
$\CO_c := \BZ_p + p^c\CO_K$.  
%If $cn = 0$, then there is an order $R$ with
%discriminant $n$ such that $R \cap K = \CO_c$. Such order is unique up to 
%$K^\times$-conjugacy. Moreover, the dimension of the invariant subspace of
%$\pi$ under $R^\times$ is less than $2$ and there is a test vector in this
%subspace. 

Assume $K$ is split or $c$ is large enough. In this case, as the functional space
$\CP(\pi,\chi)$ is nonzero, the quaternion $B$ must be split
and we shall take the test vector as the new vector of $\pi$. Precisely,
the test vector $f$ is a nonzero vector in the line of $\pi$ consisting with
invariant vectors under the actions of $R^\times$
where $R$ is an Eichler order of discriminant $n$ such that
$R \cap K = \CO_c$. By definition, all such orders are $\GL_2(\BQ_p)$-conjugate to the order $
R_0(n)$. Note that the toric integral is invariant by modifying $f$ to $\pi(t)f$ for any
$t \in K^\times$. Thus, what matters here is the $K^\times$-conjugacy classes of such Eichler orders.
In general, the set of $K^\times$-conjugacy classes of such Eichler orders is finite but not a
singleton (See also \cite{CST}, Lemma 3.2).  Here, we shall take the Eichler order $R$ to be
{\em admissible} for the pair $(\pi,\chi)$. This means that it is the intersection of two maximal orders
$R'$ and $R''$ of $M_2(\BQ_p)$ such that $R' \cap K = \CO_{c}$ and
\[R'' \cap K =
  \begin{cases}
    \CO_{c-n}, \quad &\text{if $c \geq n$};\\
    \CO_{K}, \quad &\text{otherwise.}
\end{cases}\]
Such order exists if and only if the following condition holds:
\[\tag{$*$} c-n + e - 1\geq 0 \text{ if $K$ is nonsplit}\]
where $e$ is the ramification index of $K/\BQ_p$.
If exists, it is unique up to conjugation by
normalizers of $K^\times$ in $\GL_2(\BQ_p)$. Moreover, 
it is even unique up to $K^\times$-conjugation
unless $K$ is split and $0 < c < n$. Finally, using Kirillov model of $\pi$, 
we can prove that the toric integral $\beta^0(f)$ is nonvanishing for any nonzero
$f \in \pi$ invariant under the actions of $R^\times$ with $R$ an admissible order for $(\pi,\chi)$ (See also
\cite{CST} Proposition 3.12 and Lemma \ref{lemma-toric} in this paper).

Globally, write $N = N_0N_1$ with $N_0 = \prod_{p|(N, c)} p^{\ord_p(N)}$ the $c$-part of $N$.
In particular, $(N_1,c) = 1$. Assume:
\begin{itemize}
\item The condition (1) and (2) hold for $N_1$ and $c$: 
	 if $p|(N_1,D)$ then $\ord_p(N_1) = 1$ and
	      the set  $\Sigma$  has even cardinality, where $\Sigma$ consists of prime factors $p$ of $N_1$ satisfying
		\begin{itemize}
			\item either $p$ is inert in $K$ and $\ord_p(N_1)$ is odd,
			\item or  $p$ is ramified in $K$ and $\chi([\fp]) = a_p$
				where $[\fp]\in \Pic(\CO_c)$ is the ideal class of the prime ideal $\fp$ above $p$.
		\end{itemize}
\item The condition $(*)$ holds for the $N_0$-part:  $\ord_p(c)-\ord_p(N) + e_p - 1\geq 0$ if $p|N_0$ and $K_p$ is nonsplit.
\end{itemize}

This condition implies that the root number $\epsilon(\phi,\chi)$ is $-1$ (See also Lemma \ref{epsilon}).
Let $B$ be the indefinite quoternion algebra over $\BQ$ ramified exactly at places in $\Sigma$.
Fix an embedding from $K$ to $B$. Let $R$ be an order of $B$ with discriminant $N$ such that
$R \cap K = \CO_c$ and $R_p$ is admissible for $(\phi,\chi)$ for any prime $p|N_0$.
As before, denote by
$\Gamma := \left\{ \gamma \in R^\times | \RN_B(\gamma) = 1 \right\}$ and $X_\Gamma$ is the Shimura curve
over $\BQ$ with level $\Gamma$. 
There exists an admissible modular parametrization $f_1: X_\Gamma \ra A$ which
is unique in the sense as before.
Let $P \in X_\Gamma(H_c)$ be the point represented by
the point in $\CH$ fixed by $K^\times$. Let
\[P_\chi(f_1) = \sum_{\sigma \in \Gal(H_c/K)} f_1(P)^\sigma \chi(\sigma) \in A(H_c)_\BC.\]
Similarly, there is an admissible modular parametrization $f_2: X_\Gamma \ra A^\vee$
and one may define $P_{\chi^{-1}}(f_2) \in A^\vee(H_c)_\BC$.
\begin{thm}\label{th:3}
  Assume conditions in  (1) and (2) hold for $N_1$ and $c$ while $(*)$ holds for the $N_0$-part,
  we have the following identity
  \[L'^{(S)}(1,\phi,\chi) = 2^{-\mu(N_1, D)}\cdot
    \frac{8\pi^2(\phi,\phi)_{\Gamma_0(N)}}{u^2\sqrt{|Dc^2|}}\cdot
   \frac{\langle P_\chi(f_1),P_{\chi^{-1}}(f_2)\rangle_K}{(f_1,f_2)_\Gamma}.\]
Here  $\mu(N_1, D)$ is the number of common prime factors of $N_1$ and $D$, the set 
\[S=\big\{p|(N, Dc) \big|\ \text{if $p||N$ then $\ord_p(c/N) \geq 0$}\big\}\]
and $L^{(S)}(s, \phi, \chi)$ is obtained by the L-function $L(s, \phi, \chi)$
with Euler factors at places in $S$ removed,   and other terms are the same as in Theorem \ref{gz-disjoint}.
\end{thm}

This is a special case of Theorem \ref{gz-sharp}.
The difference between the above construction of Heegner point and
the one in \cite{CST} is at those places $p|N_0$ ramified in $K$
such that $\ord_p(c) + 1 = \ord_p(N)$.
At these places, we require here that
$f$ is of $\Gamma_0(N)$-level structure while
in \cite{CST}, it is considered to be $\chi^{-1}$-eigen.
A more general construction
for abelian varieties over totally real fields is considered.
In particular, an explicit Gross-Zagier (see also Theorem \ref{gz-sharp})
is obtained as an application
of the variant of Gross-Zagier formula in \cite{CST} Theorem 1.6.

%Precisely,
%to prove $2n$ is a cube sum with $n = p^*$ and
%$p$ is a prime $\equiv 2,5 \pmod{9}$, consider
%the elliptic curve $E: y^2 = x^3 + 1$ and for any $a \in \BQ^\times$, let
%$E^{(a)}: y^2 = x^3 + a^2$. Then $E^{(a)}$ is isogenous to $C^{(a)}$ over $\BQ$.
%For the $3$-Selmer groups, we have
%$\dim_{\BF_3} \Sel_3(E^{(n)})/E^{(n)}(\BQ)_\tor \leq 1$
%while $\dim_{\BF_3} \Sel_3(E^{(n^{-1})})/E^{(n^{-1})}(\BQ)_\tor = 0$.
%Then by the BSD conjecture,
%$L'(1,E,\chi_n) = L'(1,E^{(n)})L(1,E^{(n^{-1})})$ is nonvanishing. Here,
%$K = \BQ(\sqrt{-3})$ and $\chi_n$ is the character
%taking any $\sigma \in \Gal(K^{\mathrm{ab}}/K)$ to $(\sqrt[3]{n})^{\sigma-1}$.
%In particular, $N = 36$ and $c = 3p$ with $(c,N) = 3$.

\par Next, we explain how to prove Theorem \ref{th:1} and \ref{th:2}.
Let $E/\BQ$ be the elliptic curve with Weierstrass equation $y^2=x^3+1$.
The elliptic curve $E$ is isomorphic to the modular curve $X_0(36)$ over $\BQ$.
Let $K\subset \BC$ be the imaginary quadratic field over $\BQ$ generated  by $\omega=e^{2\pi i/3}$ and $\CO_K$ its ring of integers.
It is known that $E$ is endowed with complex multiplication $[\ \ ]:\CO_K\simeq \End_\BC(E)$. For any $n\in \BQ^\times$,
let $E^{(n)}/\BQ$ be the elliptic curve $y^2=x^3+n^2$.  Then $E^{(n)}$ is isogenous to $C^{(n)}$ over $\BQ$
and $\rank_{\CO_K}\,E^{(n)}(K)=\rank_\BZ \,E^{(n)}(\BQ)$. Thus $2n$ is a cube sum if and only if $E^{(n)}(K)$ has a point of infinite order.
\par Let $k\geq 1$ be an integer. Let $p_1,p_2,\cdots,p_k$ be distinct odd primes $\equiv 2,5\mod 9$.
Let $n={p_1^{*}}^{\epsilon_1}p_2^{*\epsilon_2}\cdots p_k^{*\epsilon_k}$ with $\epsilon_i=\pm 1$ and
let $N=p_1p_2\cdots p_k$. Let $\chi_n:\Gal(K^\mathrm{ab}/K)\rightarrow \BC^\times$ be the character such
that $\chi_n (\sigma)=(\sqrt[3]{n})^{\sigma-1}, \forall \sigma\in \Gal(K^\mathrm{ab}/K)$.
Then $c(\chi) = 3N$. Denote by $E(K(\sqrt[3]{n}))^{\chi_n}$ the subgroup of $E(K(\sqrt[3]{n}))$ consisting of points $P$
such that $P^\sigma=[\chi_n(\sigma)]P$ for any $\sigma\in \Gal(K(\sqrt[3]{n})/K)$.
The group $E^{(n)}(K)$ is isomorphic to $E(K(\sqrt[3]{n}))^{\chi_n}$ under the map $\phi:E^{(n)}\rightarrow E$ given as
$(x,y)\mapsto\left (-\frac{\sqrt[3]{n}x}{n},\frac{y}{n}\right ).$ Thus $2n$ is a cube sum if and only if there is
a point in $E(K(\sqrt[3]{n}))^{\chi_n}$ with infinite order.

\par The indefinite quaternion algebra $B$ determined by $E$ and $\chi_n$ as above is $M_2(\BQ)$.
We carefully take an embedding $\rho: K \hookrightarrow B$ such that
$R_0(36) \cap K = \CO_{6N}$ with $R_0(36) =
\begin{pmatrix}
	\BZ & \BZ \\
	36\BZ & \BZ
\end{pmatrix}$
and $R_0(36)$ is admissible for $(E,\chi_n)$ at the place $3$.
A concrete such embedding is given before Theorem \ref{th:4}.
%\[\rho(\omega) =
%  \begin{pmatrix}
%    4 & -\frac{7N}{6} \\
%    \frac{18}{N} & -5
%\end{pmatrix}\]
%with fixed point $h_0 = \left(1+\frac{\sqrt{-3}}{9}\right)\frac{N}{4}$ in the upper half plane $\CH$.
%The order $\CR$ can be taken as
%\[\CR=\left \{\begin{pmatrix}
%    a & b \\
%    c & d
%\end{pmatrix}\in \M_2(\BQ)\left |\begin{array}{c}a,c/2,d,2b\in \BZ,c\equiv 0\mod 9,\\a+c/2+d\equiv 0\mod 2,\\ 2b+c/2\equiv 0\mod 2.
%\end{array} \right .\right \}\]
Let $f:X_0(36)\rightarrow E$
be the modular parametrization of degree one mapping the cusp $[\infty]$ to the identity $O$ on $E$.
Then $f$ is essentially a test vector described as above in this special situation (See also Lemma \ref{test-vector}).
\par  The L-function $L(s,E,\chi_n)$ satisfies $L(s, E, \chi_n)=L(s, E^{(n)})L(s, E^{(n^{-1})}).$
Assume $k$ is odd. Then the condition $\sum_{1}^k\epsilon_i\equiv 1\mod 3$ is equivalent to
$\epsilon(E,\chi_n)=\epsilon(E^{(n)})=-1$.
Let $\Omega^{(n)}$ denote the minimal real period of $E^{(n)}$ and $\Omega$ the one of $E$.
Let $P_0$ be the CM point on $X_0(36)$ represented by the fixed point on $\CH$ of $K^\times$
under our choice of embedding. Then $P_0$ is defined over $H_{6N}$.
Applying Theorem $\ref{th:3}$ to $E,\chi_n$ and $f$, we have the following height formula.
\begin{cor}\label{co:2}
Let $k$ be an odd integer. Suppose $p_1,p_2,\cdots,p_k$ are distinct odd prime numbers $\equiv 2,5\mod 9$.
Let $n={p_1^{*}}^{\epsilon_1}p_2^{*\epsilon_2}\cdots p_k^{*\epsilon_k}$ with $\epsilon_i=\pm 1$ and
$\sum \epsilon_i \equiv 1 \mod 3$. The Heegner divisor
$$z_{n}:=\sum_{\sigma\in \Gal(H_{6N}/K)}[\chi_n^{-1}(\sigma)]f(P_0)^{\sigma} \in E(K(\sqrt[3]{n}))^{\chi_n}$$
satisfies the following Gross-Zagier formula:
\[L'(1,E^{(n)})L(1, E^{(n^{-1})})=3^{-3}\cdot \Omega^{(n)}\Omega^{(n^{-1})}\cdot \widehat{h}_\BQ(z_{n}),\]
where $\wh{h}_\BQ$ is the canonical \Neron-Tate height on $E(\bar{\BQ})$ over $\BQ$.
\end{cor}
From this height formula, we see that the Heegner point $z_n$ is of infinite order if and only if the L-function $L(s,E,\chi_n)$ has vanishing order $1$ at $s=1$. Recall when $k=1$, $\dim_{\BF_3}\Sel_3(E^{(p^*)})/E^{(p^*)}(\BQ)_\tor\leq 1$ and $\dim_{\BF_3}\Sel_3(E^{(p^{*-1})})/E^{(p^{*-1})}(\BQ)_\tor=0$. The Birch and Swinnerton-Dyer conjecture predicts that $$\ord_{s=1} L(s,E^{(p^*)})=\rank_\BZ\,E^{(p^*)}(\BQ)=1 \textrm{ and } \ord_{s=1} L(s,E^{(p^{*-1})})=\rank_\BZ\,E^{(p^{*-1})}(\BQ)=0.$$
We will prove this by showing that the Heegner point $z_{p^*}$ is of infinite order.
\par Next we explain how to prove the nontriviality of Heegner points in general. For any  $d=p_1^{\epsilon_1}p_2^{\epsilon_2}\cdots p_k^{\epsilon_k}$ with $\epsilon_i=0,\pm 1$,  let $\chi_d:\Gal(H_{6N}/K)\rightarrow \BC^\times$ be the character such that $\chi_d (\sigma)=(\sqrt[3]{d})^{\sigma-1}, \forall \sigma\in \Gal(H_{6N}/K)$.  Define the Heegner divisor
$$z_{d}:=\sum_{\sigma\in \Gal(H_{6N}/K)}[\chi_d^{-1}(\sigma)]f(P_0)^{\sigma}\in E(K(\sqrt[3]{d}))^{\chi_d}.$$

\par Summing over all these Heegner divisors, we have an identity:
$$\sum_{d}z_{d}=3^kz_0,$$ where $z_0=\Tr_{H_{6N}/H_0}f(P_0)\in E(H_0)$ is the genus Heegner point and $H_0=K(\sqrt[3]{p_1},\sqrt[3]{p_2},\cdots,\sqrt[3]{p_k})\subset H_{3N}$ is the genus field over $K$.
By a similar argument due to Heegner (see also Birch as in \cite{Birch}, and Satge in \cite{Satge}), it follows from the key fact $3\nmid [H_{3N}: H_0]$ that the genus Heegner point $z_0$ is of infinite order. Then we examine which Heeger points will contribute on the left hand side of the equality. By Euler system property  of Heegner divisors, for any $d=p_1^{\epsilon_1}p_2^{\epsilon_2}\cdots p_k^{\epsilon_k}$ with some $\epsilon_i=0$ , the Heegner divisor $z_{d}$ is zero. Thus only Heegner points $z_d$ with $\epsilon_i\neq 0,\forall i$, remain and at least one such Heegner point is of infinite order and Theorem $\ref{th:1}$ follows. Note that if $k=1$, there is exactly one term remains on the left hand side and the first part of Theorem $\ref{th:2}$ follows. Comparing the height formula in Corollary $\ref{co:2}$ and BSD conjecture, the $3$-part of the second assertion in Theorem $\ref{th:2}$ follows from the $3$-divisibility of the Heegner point $z_{n}$ in $E(K(\sqrt[3]{n}))^{\chi_n}$. The $\ell$-part of the second assertion in Theorem $\ref{th:2}$, $\ell\nmid 6p$, follows from work of Kolyvagin \cite{K}, Perrin-Riou \cite{Perrin-Riou} and Kobayashi \cite{Kobayashi}

\section{Nontriviality of Heegner Points}
Recall the elliptic curve $E/\BQ$ is given by the Weierstrass equation $y^2=x^3+1$ and
the complex multiplication is fixed by the isomorphism $[\,\,\,]:\CO_K\rightarrow \End_\BC(E)$ such that $[\omega](x,y)=(\omega x,y)$.
\begin{prop}\label{pro:1}
The group of rational points $E(\BQ)$ is torsion cyclic of order $6$
and is generated by the point $(2,3)$. The group of $K$-points
$E(K)$ is torsion and equal to $E[2\sqrt{-3}]$.
\end{prop}
\begin{proof}
Let $m$ be an integer such that $[m]E(\BQ)_{\mathrm{tor}}=0$. For
any prime number $p\nmid 6m$, $E(\BQ)_{\tor}\subset
E(\BQ_p)[m]\simeq E(\BF_p)[m]$ under reduction $\mod p$. When
$p\equiv 5\mod 6$, the reduction of $E \mod p$ is supersingular, and
then $E(\BF_p)=p+1$. Thus $\#E(\BQ)_\tor|6$. On the other hand, the
point $(2,3)$ is a rational point of order $6$. So $E(\BQ)_\tor$ is
cyclic of order $6$ and generated by $(2,3)$. There is a rational
$2$-torsion point $(-1,0)$ on $E$ and, using the $2$-descent method,
one can show $\dim_{\BF_2}\mathrm{Sel}_2(E/\BQ)/E[2](\BQ)=0$ and
hence $\rank_\BZ \,E(\BQ)=0$.
\par For any point $P\in E(K)$, $P+\bar{P}$ is a rational point and $[6](P+\bar{P})=0$. If $[6]P\neq O$, then it has coordinates of form $(x,y\sqrt{-3})$ with $x,y\in \BQ$ and so does $[6\omega]P=(\omega x,y\sqrt{-3})$. But this forces $x=0$ which is impossible. So $E(K)$ is torsion.
\par The conductor of the Hecke character associated to the elliptic curve $E/K$ is the ideal $(2\sqrt{-3})$. By the main theorem of the complex multiplication, $E(K)\subset E[2\sqrt{-3}]$. On the other hand, there are exactly $12$ $K$-points generated by $E(\BQ)$ under complex multiplication by $\CO_K$, namely, $O,(-\omega^n,0)$, $(0,\pm 1)$, $(2\omega^n,3),$ $(2\omega^n,-3)$ with $n=0,1,2$. So $E(K)_{\mathrm{tor}}=E[2\sqrt{-3}]$.
\end{proof}
\par Let $U_0(36)$ be the open subgroup of $\GL_2(\widehat{\BZ})$ consisting of matrices
$\left (
\begin{array}{cc}
a&b\\c&d
\end{array}
\right )$ such that $c\equiv 0\mod 36$. Let
$\Gamma_0(36)=\GL_2(\BQ)^+\cap U_0(36)$. Then $X_0(36)$ is the
modular curve over $\BQ$ whose underlying Riemann surface is
\[X_0(36)(\BC)=\GL_2(\BQ)^+\backslash\left (\CH\sqcup\BP^1(\BQ)\right )\times \GL_2(\BA_f)/U_0(36)\simeq \Gamma_0(36)\backslash \CH\bigsqcup \Gamma_0(36)\backslash \BP^1(\BQ).\]
Elements in $S=\Gamma_0(36)\backslash \BP^1(\BQ)$ are called cusps
of $X_0(36)$ and there are $12$ cusps, $6$ among which are rational
and the other $6$ are defined over the imaginary quadratic field
$K=\BQ(\sqrt{-3})$. For each $z\in \BP^1(\BQ)$, denote by $[z]$ the
cusp on $X_0(36)$ represented by $z$. The modular curve $X_0(36)$ is
a projective smooth curve of genus $1$ over $\BQ$ with a rational
cusp $[\infty]$ and we have an elliptic curve $(X_0(36),[\infty])$
over $\BQ$ with $[\infty]$ as its zero element.
\par For any field extension $F/\BQ$, we write $\Isom_F(X_0(36))$ for the group of algebraic isomorphisms of $X_0(36)$ defined over $F$, and $\Aut_F(X_0(36))$ for the subgroup of algebraic isomorphisms over $F$ which fix the cusp $[\infty]$. Define $N$ to be the normalizer of $\Gamma_0(36)$ in $\GL_2^+(\BQ)$.
The action of $N$ on $X_0(36)$ induces an injective homomorphism
\[T:N/\BQ^\times\Gamma_0(36)\rightarrow \Isom_\BC(X_0(36))=\Aut_\BC(X_0(36))\ltimes X_0(36)(\BC).\]
\begin{prop}\label{pro:2}
The elliptic curve $(X_0(36),[\infty])$ has complex multiplication
by $\CO_K$ and Weierstrass equation $y^2=x^3+1$ such that the cusp
$[0]$ has coordinates $(2,3)$.
\end{prop}
We identify the elliptic curve $(X_0(36),[\infty])$ with $E$ and
then $E(K)=S$.
\begin{proof}There are two special linear fractional transformations $X_0(36)$ given by matrices
\[A=\left (\begin{array}{cc}1&1/6\\0&1\end{array}\right )\textrm{ and }B=\left (\begin{array}{cc}0&1\\-36&0\end{array}\right ).\]
Since the transformation $T(A)$ is of order $6$ and fixes
$[\infty]$, $T(A)$ acts on the elliptic curve $(X_0(36),[\infty])$
as a complex multiplication by a primitive sixth root of unity. Thus
the elliptic curve $(X_0(36),[\infty])$ has complex multiplication
by $\CO_K$. The linear fractional transformation $T(B)$ is the
Atkin-Lehner involution on $X_0(36)$, which is rational of order
$2$.
\par For $\epsilon\in \Aut_\BC(X_0(36))=\CO_K^\times$ and $\alpha\in X_0(36)(\BC)$, denote $t_{\epsilon,\alpha}$ to be the isomorphism $P\mapsto [\epsilon]P+\alpha$ in $\Aut_\BC(X_0(36))\ltimes X_0(36)(\BC)$. We simply write $t_{1,\alpha}$ as $t_\alpha$, and $t_{\epsilon, [\infty]}$ as $[\epsilon]$. Note $t_{\epsilon,\alpha}^2=t_{\epsilon^2,[\epsilon]\alpha+\alpha}$ for any $\epsilon\in \CO_K^\times$ and $\alpha\in X_0(36)(\BC)$. Since $T(B)$ is rational of order $2$, it is either $t_{-1,\alpha}$ for some rational $\alpha$, or $t_{\alpha}$ for some rational point $\alpha$ of order $2$. Note $T(B)$ fixes the point $[\sqrt{-1}/6]$ and then $T(B)=t_{-1,\alpha}$ for some rational $\alpha$. Since $t_{\alpha}=t_{-1,\alpha}\circ [-1]=T(BA^3)$ and the matrix
\[BA^3=\left (\begin{array}{cc}
0&1\\-36&-18
\end{array}\right )\]
is of order $6$ in $\Phi$ and takes cusp $[\infty]$ to cusp $[0]$,
the cusp $[0]$ is rational of order $6$ and the Atkin-Lehner
involution $T(B)=t_{-1,[0]}$.
\par It is known that every elliptic curve over $\BQ$ is parametrized by the modular curve of the same level as its conductor. It follows that $(X_0(36),[\infty])$ is isogenous to the elliptic curve $E:y^2=x^3+1$. However, there are exactly $4$ isomorphism classes of elliptic curves over $\BQ$ with conductor $36$ and the elliptic curve $y^2=x^3+1$ is the unique one with a rational point of order $6$ and with complex multiplication by $\CO_K$. So $(X_0(36),[\infty])$ is isomorphic to the elliptic curve $E:y^2=x^3+1$. The rational points of order $6$ on $y^2=x^3+1$ are $(2,\pm 3)$. This isomorphism is unique if it maps $[0]$ to $(2,3)$. Thus $E$ has a Weierstrass equation $y^2=x^3+1$ with $(2,3)$ as the coordinates of the cusp $[0]$.
\end{proof}
\par Let $\phi=\sum_{n\geq 1} a_nq^n$ be the unique new form of level $\Gamma_0(36)$ and weight $2$. Then the N\'eron differential on $X_0(36)$ is $dx/2y=\phi(q)dq/q$. At $[\infty]$, it is represented by $dq$ and $T(A)^*(dq)=-\omega^2dq$. On the other hand, $[-\omega^2](x,y)=(\omega^2,-y)$ and then $[-\omega^2]^*(dx/2y)=-\omega^2 dx/2y$. So we have $T(A)=[-\omega^2]$.
\begin{prop}\label{pro:3}
The natural homomorphism $T$ induces an isomorphism $$N/\BQ^\times\Gamma_0(36)\simeq
\Isom_K(E)= \Aut_K(E)\ltimes E(K).$$ More precisely, one has the
following relations:
\[\begin{aligned}
&t_O=T\left (\begin{array}{cc}1&0\\0&1\end{array}\right ),&
&t_{(0,1)}=T\left (\begin{array}{cc}-2&-1\\36&16\end{array}\right
),&
&t_{(0,-1)}=T\left (\begin{array}{cc}16&1\\-36&-2\end{array}\right ),&\\
&t_{(-1,0)}=T\left (\begin{array}{cc}9&4\\-144&-63\end{array}\right
),& &t_{(-\omega,0)}=T\left
(\begin{array}{cc}39&4\\144&15\end{array}\right ),&
&t_{(-\omega^2,0)}=T\left (\begin{array}{cc}87&-20\\144&-33\end{array}\right ),&\\
&t_{(2,3)}=T\left (\begin{array}{cc}0&1\\-36&-18\end{array}\right
),& &t_{(2\omega,3)}=T\left
(\begin{array}{cc}12&1\\36&6\end{array}\right ),&
&t_{(2\omega^2,3)}=T\left (\begin{array}{cc}12&11\\-36&-30\end{array}\right ),&\\
&t_{(2,-3)}=T\left (\begin{array}{cc}-18&-1\\36&0\end{array}\right
),& &t_{(2\omega,-3)}=T\left
(\begin{array}{cc}-6&1\\36&-12\end{array}\right ),&
&t_{(2\omega^2,-3)}=T\left
(\begin{array}{cc}6&1\\36&12\end{array}\right ).&
\end{aligned}
\]
\end{prop}
\begin{proof}
For any $C\in \GL_2(\BQ)^+$, if $T(C)[\infty]=\alpha\in S=E(K)$,
then $t_{-\alpha}\circ T(C)$ is an element in $\Aut_K(E)$. Thus the
image of $T$ lies in $\Aut_K(E)\ltimes E(K)$. Since
$T(A)=[-\omega^2]$, $T(BA^3)=t_{[0]}$ and $[-\omega^2], t_{[0]}$
generates $\Isom_K(E)$, the homomorphism $T$ is an isomorphism
between $N/\BQ^\times\Gamma_0(36)$ and $\Isom_K(E)$.
\par The verification of the relations is straightforward.
\end{proof}
\begin{cor}
The complete list of cusps on $X_0(36)$ is given as follows:
\[[0],[1/2],[1/3],[-1/3],[-1/16],[1/6],[-1/6],[-4/9],[13/48],[29/48],[-1/18],[\infty].\]
\end{cor}
As is seen, we have $S\simeq E[2\sqrt{-3}]$. Let
\[\tau:\BZ[\omega]/(2\sqrt{-3})\longrightarrow S\]
be the unique isomorphism such that $\tau(0)=[\infty]$ and
$\tau(1)=[0]$. Then we have that
\[\begin{aligned}
&\tau(-1)=[-1/2],\quad &\tau(\omega )=[1/3],\quad
&\tau(\omega^2)=[-1/3],\quad &\tau(3)=[-1/16],\\
&\tau(-\omega)=[-1/6],\quad
&\tau(-\omega^2)=[1/6],\quad&\tau(4)=[-4/9],\quad &
\tau(3\omega)=[13/48],\\&\tau(3\omega^2)=[29/48],\quad
&\tau(2)=[-1/18].\quad&
\end{aligned}\]
\par Let $U\subset U_0(36)$ be the subgroup of index $2$ consisting of matrices
$$\left (\begin{array}{cc}a&b\\c&d
\end{array} \right ) \textrm{ with $a\equiv d\mod 3$.}$$
Let $X_U$ be the modular curve over $\BQ$ whose underlying Riemann
surface is
\[X_U(\BC)=\GL_2(\BQ)^+\backslash\left (\CH\sqcup\BP^1(\BQ)\right )\times \GL_2(\BA_f)/U.\]
Under class field theory, $ \BA_f^\times/\BQ^\times\det(U)\simeq \Gal(K/\BQ).$
By noting that $\GL_2(\BQ)^+\cap U=\Gamma_0(36)$, we see that the
modular curve $X_U$ is isomorphic to $X_0(36)\times_\BQ K$ as a curve over $\BQ$. Let
$U_0(36)/U=\langle\epsilon\rangle$ where
\[\epsilon=\left (\begin{array}{cc}1&0\\0&-1\end{array}\right ).\]
The non-trivial Galois action of $\Gal(K/\BQ)$ on $X_U$ is given by
the right-translation of $\epsilon$ on $X_U$. We have
\[\Isom_\BQ(X_U)=\left (\BZ[\omega]/(2\sqrt{-3})\rtimes\CO_K^\times\right )\rtimes\Gal(K/\BQ).\]
Let $N_{\GL_2(\BA_f)}(U)$ be the normalizer of $U$ in
$\GL_2(\BA_f)$. Then there is a natural homomorphism
$N_{\GL_2(\BA_f)}(U)\longrightarrow \Isom_\BQ(X_U)$ induced by right
translation on $X_U$. The curve $X_U$ is not geometrically connected
and has two connected components over $\BC$. An element $g\in
N_{\GL_2(\BA_f)}(U)$ maps one component of $X_U$ onto another if and
only if it has image $-1$ under the composition of the following
morphisms:
\[\xymatrix{\GL_2(\BA_f)=\GL_2(\BQ)^+U_0(32)\ar[r]^{\quad\quad \quad\quad\quad \det}&\BQ^{\times}_+\widehat{\BZ}^\times\ar[r]&\BZ_3^\times/(1+3\BZ_3),}\]
where the last morphism is the projection from
$\widehat{\BZ}^\times$ to its $3$-adic factor.

Let $k\geq 1$ be an odd integer. Let $p_1, p_2,\cdots,p_k$ be distinct odd primes $\equiv 2, 5\mod
9$ and $N=p_1p_2\cdots p_k$. Let $\rho: K \ra M_2(\BQ)$ be the normalized embedding given by
\[\omega \mapsto
  \begin{pmatrix}
    4 & -\frac{7N}{6} \\
    \frac{18}{N} & -5
\end{pmatrix}\]
with fixed point $h_0 = \left(1+\frac{\sqrt{-3}}{9}\right)\frac{N}{4} \in \CH$. Then
$\widehat{K}^\times \cap U = \widehat{\CO}_{6N}^\times$. Let $w=\begin{pmatrix}
    1 & -\frac{N}{2} \\
    0 & -1
\end{pmatrix}$ be an element in $N_{\mathrm{GL}_2(\BQ)^+}(K^\times)$ which normalises $U$ and hence induces an isomorphism in $\Isom_\BQ(X_U)$. The groups $K_2^\times/\BQ_2^\times(1+2\CO_2)=\langle\omega_2\rangle^{\BZ/3\BZ}$ and $\CO_3^\times/\BZ_3^\times(1+3\CO_2)=\langle\omega_3\rangle^{\BZ/3\BZ}$. By noting that $\omega_2$ and $\omega_3$ normalises $U$ in $\GL_2(\BA_f)$, we have natural homomorphisms
\[K_2^\times/\BQ_2^\times(1+2\CO_2)\rightarrow\Isom_\BQ(X_U)\textrm{ and }\CO_3^\times/\BZ_3^\times(1+3\CO_3)\rightarrow\Isom_\BQ(X_U).\]

\begin{thm}\label{th:4}
For any point $P\in X_U$,
\[P^{\omega_2}=P+\tau(2),\,P^{\omega_3}=[\omega]P,\,P^{w\epsilon}=[-1]P.\]
\end{thm}
\begin{proof}
Since $\det\omega_2=1,\det\omega_3=1$ and $\det w\epsilon=1$, there exist $\alpha,\beta,\gamma\in \CO_K^\times$ and $R,S,T\in E(K)$ such that, for all $P\in X_U$,
\[P^{\omega_2}=[\alpha]P+R,\,P^{\omega_3}=[\beta]P+S,\textrm{ and }P^{w\epsilon}=[\gamma] P+T.\]
Note
\[
 \begin{pmatrix}
    1 & 2^{-1} \\
    18 & 10
  \end{pmatrix}\omega_2 \in U \quad \text{and} \quad
  \begin{pmatrix}
    1 & 3^{-1} \\
    36 & 13
\end{pmatrix}\omega_3 \in U, \begin{pmatrix}
1&\frac{N}{2}\\0&1
\end{pmatrix}w\epsilon\in U.\]
Taking $P=[\infty]$,
\[R=[\infty]^{\omega_2}=[1/18],S=[\infty]^{\omega_3}=[1/36],T=[\infty]^{w\epsilon}=[\infty].\]
Taking $P=[0]$,
\[[\alpha][0]+[1/18]=[1/20],[\beta][0]=[1/39],[\gamma][0]=[N/2].\] When identified with identities in $\CO_K/(2\sqrt{-3})$, it follows that $\alpha=1,\beta=\omega$ and $\gamma=-1$.
\end{proof}
\par Let $P_0=[h_0,1]$ be a complex multiplication point in $X_0(36)$. By the main theorem of complex multiplication, the point $P_0$ is in $E(K^{\mathrm{ab}})$. Let $\sigma:K^\times\backslash \widehat{K}^\times\rightarrow G_K^{\mathrm{ab}}$ be the Artin reciprocity law.
\begin{cor} \label{co:1}
The point $P_0\in E(H_{6N})$ satisfies
\[P_0^{\sigma_{\omega_2}}=P_0+\tau(2),\,P_0^{\sigma_{\omega_3}}=[\omega]P_0,\,\overline{P_0}=[-1]P_0.\]
\end{cor}
\begin{proof}
Consider the normalised embedding $\rho:K\longrightarrow M_2(\BQ)$. Let $P=[h_0,1]$ be a point in $X_U$. By Shimura's reciprocity law, $P^{\sigma_t}=[h_0,\rho(t)]$ for any $t\in \widehat{K}^\times$. Thus $P^{\sigma_t}=P$ if and only if $t\in K^\times \left (U\cap \widehat{K}^\times\right )=K^\times \widehat{\CO}_{6N}^\times$. Hence $P\in X_U(H_{6N})$. In view of Theorem $\ref{th:4}$, by Shimura's reciprocity law, one has
\[P^{\sigma_{\omega_2}}=P+\tau(2),\,P^{\sigma_{\omega_3}}=[\omega]P,\,\overline{P}=[-1]P.\]
Since $P_0$ is the image of $P$ under the natural projection, the corollary follows.
\end{proof}

Let $H_0=K(\sqrt[3]{p_1},\sqrt[3]{p_2},\cdots,\sqrt[3]{p_k})$ be the genus field over $K$.
\begin{prop}\label{pro:4}
\begin{itemize}
\item[a.] The field $H_0$ is contained in $H_{3N}$ and the Galois group $\Gal(H_0/K)=\langle\sigma_{\omega_p} \rangle_{p\mid N}$ and $\sigma_{\omega_3}(\sqrt[3]{p^*})=\omega\sqrt[3]{p^*}$ for $p\mid N$.
\item[b.] $H_6=K(\sqrt[3]{2})$ with Galois group $\Gal(H_6/K)=\langle\sigma_{\omega_2}\rangle$ and $H_{6N}=H_{3N}(\sqrt[3]{2})$ with $\Gal(H_{6N}/H_{3N})=\langle\sigma_{\omega_2}\rangle$ and
$\sigma_{\omega_2}(\sqrt[3]{2})=\sigma^{-1}_{\omega_3}(\sqrt[3]{2})=\omega^2\sqrt[3]{2}$.
\end{itemize}
\end{prop}
\begin{proof}
For $p\mid N$, let $L=K(\sqrt[3]{p})$ and $H_0$ is the composition
filed of these cubic extensions $L$. The field extension $L/K$ is
totally ramified at the prime ideals $(\sqrt{-3})$ and $(p)$. So the
local norm group of $L$ at $(\sqrt{-3})$ (resp. $(p)$) cannot
contains $\CO_3^\times$ (resp. $\CO_p^\times$). We claim that
$\BZ_3^\times(1+3\CO_3)\BZ_p^\times(1+p\CO_p)$ fixes $L$ under the
Artin reciprocity law
\par First we show that $\BZ_p^\times(1+p\CO_p)$ fixes $L$. Let $e(X)=-pX+X^{p^2-1}$ be the Lubin-Tate series of $K_p$ w.r.t. the primitive element $-p$ and let $\CF$ be the Lubin-Tate $\CO_p$-module associated to $e(X)$. For any element $\alpha\in \CO_p$, denote the action of $\alpha$ on $\CF$ as $[\alpha]_{\CF}$. Then it is easily checked that $[\omega_p]_\CF(X)=\omega X$. Let $K_p^1$ be the first Lubin-Tate extension which is the splitting field of $e(X)$ over $K_p$ and contains $\sqrt[3]{p}$. Let $u\in \BZ_p^\times (1+p\CO_p)$. Since $[u]_\CF(e(X))=e([u]_\CF(X))$, $[u]_\CF(\sqrt[p^2-1]{p})=\omega_{p^2-1}^l \sqrt[p^2-1]{p}$ for some integer $l$, where $\omega_{p^2-1}$ is a primitive $(p^2-1)$-th root of unity. If $u\equiv \omega_{p^2-1}^{(p+1)m}\mod p$ for some integer $m$, $[u]_\CF(\sqrt[p^2-1]{p})\equiv \omega_{p^2-1}^{(p+1)m}\sqrt[p^2-1]{p}\mod (\sqrt[p^2-1]{p})^2$. Since $\omega_{p^2-1}^l-\omega_{p^2-1}^{(p+1)m}$ is either zero or a $p$-adic unit, we conclude that $\omega_{p^2-1}^l=\omega_{p^2-1}^{(p+1)m}$ and $[u]_\CF(\sqrt[p^2-1]{p})=\omega_{p^2-1}^{(p+1)m}\sqrt[p^2-1]{p}$. Hence $$[u]_\CF(\sqrt[3]{p})=(\omega_{p^2-1}^{(p+1)m}\sqrt[p^2-1]{p})^{\frac{p^2-1}{3}}=\sqrt[3]{p}.$$ Thus $\BZ_p^\times(1+p\CO_p)$ fixes $L$.
\par Next we show that $N_{L_3/K_3}(\CO_{L,3}^\times)=\BZ_3^\times(1+3\CO_3)$. Let $L_3$ be the local field of $L$ at the unique prime ideal above $3$ and $\CO_{L,3}$ is the ring of integers of $L_3$. Note in $L_3$, $(1+\sqrt[3]{p})(1+\omega\sqrt[3]{p})(1+\omega^2\sqrt[3]{p})=1+p$. Since $p\equiv 2\mod 3$, the $3$-adic valuation $\nu_3(1+\sqrt[3]{p})=\frac{1}{3}$. Then $\varpi=\frac{\sqrt{-3}}{1+\sqrt[3]{p}}$ has $3$-adic valuation $\frac{1}{6}$. So the $6$-degree extension $L_3/\BQ_3$ is totally ramified and $\varpi$ is a uniformizer. The ring of integers $\CO_{L,3}=\BZ_3[\varpi]$. Let $x=\alpha+\beta\varpi+\gamma\varpi^2\in \CO_{L,3}^\times$. Then $\alpha\in \BZ_3[\sqrt{-3}]^\times$ and $\beta$, $\gamma\in \BZ_3[\sqrt{-3}]$. If $p\equiv 2\mod 9$, then
\begin{eqnarray*}
N_{L_3/K_3}(x)&\equiv&(\alpha+\gamma+\frac{\sqrt{-3}}{3}\beta)^3+(-\frac{\sqrt{-3}}{3}\beta)^3p+(-\gamma+\frac{\sqrt{-3}}{3}\beta)^3p^2\\
&&-3p(-\frac{\sqrt{-3}}{3})\beta(\alpha+\gamma+\frac{\sqrt{-3}}{3}\beta)(-\gamma+\frac{\sqrt{-3}}{3}\beta)\\
&\equiv& \sqrt{-3}\beta(\alpha^2-\beta^2) + A\mod 3\CO_3,
\end{eqnarray*}
where $A\in \BZ_3$. If $p\equiv 5\mod 9$, then
\begin{eqnarray*}
N_{L_3/K_3}(x)&\equiv&(\alpha-\frac{\sqrt{-3}}{3}\beta)^3+(\frac{\sqrt{-3}}{3}\beta-\gamma)^3p+(-\frac{\sqrt{-3}}{3}\beta-\gamma)^3p^2\\
&&+3p(\alpha-\frac{\sqrt{-3}}{3}\beta)(\frac{\sqrt{-3}}{3}\beta-\gamma)(\frac{\sqrt{-3}}{3}\beta+\gamma)\\
&\equiv& \sqrt{-3}\beta(\beta^2-\alpha^2) + B\mod 3\CO_3,
\end{eqnarray*}
where $B\in \BZ_3$. Since $\alpha\in \BZ_3[\sqrt{-3}]^\times$,
$\beta(\alpha^2-\beta^2)\in \sqrt{-3}\CO_3$. Hence
$N_{L_3/K_3}(\CO_{L,3}^\times)\subset \BZ_3^\times(1+3\CO_3)$. By
local class field theory, $\Gal(L_3/K_3)\simeq
\CO_3^\times/N_{L_3/K_3}(\CO_{L,3}^\times)$, because $(\sqrt{-3})$
is totally ramified in $L_3/K_3$. Since $\BZ_3^\times(1+3\CO_3)$ has
index $3$ in $\CO_3^\times$, we see that
$N_{L_3/K_3}(\CO_{L,3}^\times)=\BZ_3^\times(1+3\CO_3)$.
\par Since $L/K$ is unramified outside $(\sqrt{-3})$ and $(p)$, $K^\times \BZ_3^\times(1+3\CO_3) \CO_p^{\times 3}(1+p\CO_p)\widehat{\CO}_K^{\times (3p)}$ is contained in $K^\times N_{L/K}(\widehat{L}^\times)$. Since it has index $3$ in $\widehat{K}^\times$, $K^\times N_{L/K}(\widehat{L}^\times)=K^\times \BZ_3^\times(1+3\CO_3)\CO_p^{\times 3}(1+p\CO_p)\widehat{\CO}_K^{\times (3p)}$ and the Galois group $\Gal(L/K)$ is isomorphic to $\CO_3^\times/\BZ_3^\times(1+3\CO_3)=\langle\sigma_{\omega_3}\rangle$ and is also isomorphic to $\CO_p^\times/\CO_p^{\times 3}(1+p\CO_p)=\langle\sigma_{\omega_p}\rangle$ with $\sigma_{\omega_3}\sigma_{\omega_p}=1$ on $L$.
The Artin symbol
$[\omega_p,K_p^1/K_p](\sqrt[3]{p})=[\omega_p^{-1}]_\CF(\sqrt[3]{p})=(\omega^{-1}\sqrt[p^2-1]{p})^{\frac{p^2-1}{3}}.$
If $p\equiv 2\mod 9$,
$\sigma_{\omega_p}(\sqrt[3]{p})=\omega^2\sqrt[3]{p}$. If $p\equiv
5\mod 9$, $\sigma_{\omega_p}(\sqrt[3]{p})=\omega\sqrt[3]{p}$. Then
$\sigma_{\omega_3}(\sqrt[3]{p^*})=\sigma_{\omega_p}^{-1}(\sqrt[3]{p^*})=\omega(\sqrt[3]{p^*})$.
Since $p\equiv 2\mod 3$, $\BZ_p^\times\subset \CO_p^{\times
3}(1+p\CO_p)$ and then $\widehat{\CO}^\times_{3N}$ fixes
$\sqrt[3]{p}$. So $\sqrt[3]{p}$ is in $H_{3p}$. Thus $H_0\subset
H_{3N}$ and the norm subgroup of $H_0$ is $K^\times
\BZ_3^\times(1+3\CO_3)\prod_{p\mid
N}\CO_p^{\times 3}(1+p\CO_p)\widehat{\CO}_K^{\times (3N)}$, and hence,
the Galois group $\Gal(H_0/K)=\langle\sigma_{\omega_p} |p\mid
N\rangle$.
\par The field extension $K(\sqrt[3]{2})$ is totally ramified at $2$ and $3$ and unramified elsewhere. It can be showed in a completely same manner with the prime $p$ replaced by $2$ that $\widehat{\CO}_{6}^\times=(1+2\CO_2)\BZ^\times_3(1+3\CO_3)\widehat{\CO}_K^{\times (6)}$ fixes $\sqrt[3]{2}$. Since $K^\times \widehat{\CO}_{6}^\times$ has index $3$ in $\widehat{K}^\times$, $K(\sqrt[3]{2})$ is the class field $H_6$ and $\Gal(H_6/K)\simeq\CO_2^\times/(1+2\CO_2)=\langle\sigma_{\omega_2}\rangle$ and $\sigma_{\omega_2}(\sqrt[3]{2})=\omega^2\sqrt[3]{2}$. Since the field extension $H_{3N}/K$ is unramified at the prime ideal $(2)$, $\sqrt[3]{2}$ cannot be contained in $H_{3N}$ and, hence $H_{6N}=H_{3N}(\sqrt[3]{2})$ and $\Gal(H_{6N}/H_{3N})\simeq \CO_2^\times/(1+2\CO_2)=\langle\sigma_{\omega_2}\rangle$.
\end{proof}

Let $T=(-\sqrt[3]{4},\sqrt{-3})$ be a point in $E[3](H_6)$. By Proposition
$\ref{pro:4}$,
\[T^{\sigma_{\omega_2}}=T^{\sigma^{-1}_{\omega_3}}=[\omega]T=T+\tau(2).\]
Then $P_0-T$ is a point in $E(H_{3N})$ by Proposition
$\ref{pro:4}$ and Corollary $\ref{co:1}$.
Let $y_0=\RTr_{H_{3N}/H_0}(P_0-T)\in E(H_0)$.
\begin{prop} \label{pro:5}
\begin{itemize}
\item[a.] The point $y_0\in E(H_0)$ satisfies
\[y_0^{\sigma_{\omega_3}}=[\omega] y_0+t,\]
where $t\in E(\BQ)[3]$ is a nonzero point, and
\[\overline{y_0}= [-1]y_0.\]
\item[b.] The point $y_0\in E(H_0)$ is of infinite order. In particular, the point
\[z_0=\RTr_{H_{6N}/H_0}P_0=3y_0\]
is of infinite order in $E(H_0)$.
\end{itemize}
\end{prop}
\begin{proof}
We have
\begin{eqnarray*}
y_0^{\sigma_{\omega_3}}&=&\RTr_{H_{3N}/H_0}(P_0^{\sigma_{\omega_3}}-T^{\sigma_{\omega_3}})\\
&=&\RTr_{H_{3N}/H_0}\left [([\omega]P_0)-[\omega]^2T\right ]\\
&=&\RTr_{H_{3N}/H_0}\left [[\omega](P_0-T)+[\omega-\omega^2]T\right ]\\
&=&\RTr_{H_{3N}/H_0}\left [[\omega](P_0-T)+\tau(4)\right ]\\
&=&[\omega]y_0+t,
\end{eqnarray*}
where $t=\prod_{p\mid N}\frac{p+1}{3}\tau(4)$ is a torsion point of order $3$, because $\prod_{p\mid N}\frac{p+1}{3}$ is prime to $3$. Clearly, one has  $\overline{y_0}=[-1]y_0$.
The first assertions follows.
\par Next we prove the second assertion. We show that $E(H_0)[3^\infty]=E(K)[3^\infty]$. Assume $E(H_0)[3^\infty]\neq E(K)[3^\infty]$ and let $P$ be a point in $E(H_0)[3^\infty]\backslash E(K)$. Then $K(P)/K$ is unramified outside $2,3$. Meanwhile, $K(P)$ must contains some $\sqrt[3]{d}$, with some $d=p_1^{\epsilon_1}p_2^{\epsilon_2}\cdots p_k^{\epsilon_k}\neq 1$. Then $K(P)/K$ is ramified at the primes $p_i$ with $\epsilon_i\neq 0$. This is a contradiction. So $E(H_0)[3^\infty]=E(K)[3^\infty]$.
\par Let $C$ be a sufficiently large integer prime to $3$ that kills points in $E(H_0)_\tor$ of order prime to $3$. If $y_0$ is torsion, then $Cy_0\in E(H_0)[3^\infty]=E(K)[3^\infty]$ satisfies
\[Cy_0=Cy_0^{\sigma_{\omega_3}}=[\omega]Cy_0+ Ct,\]
where $Ct\in E(\BQ)[3]$ is nonzero by the first assertion. Then the torsion point $t$ is divisible by $[1-\omega]$ which implies that it is killed by $[2]$ because $E(K)=E[2\sqrt{-3}]$ and this conflicts with the fact the point $t$ has order $3$.
Hence, $y_0$ is a point of infinite order in $E(H_0)$ as required.
\par Finally, $z_0=\RTr_{H_{6N}/H_0}(P_0-T+T)=3y_0+\RTr_{H_{6N/H_0}}T$. Since $\Gal(H_{6N}/H_{3N})=\langle \sigma_{\omega_2}\rangle$, \[\RTr_{H_{6N}/H_{3N}}T=T+T^{\sigma_{\omega_2}}+T^{\sigma_{\omega_2}^2}=[1+\omega+\omega^2]T=0.\]
Hence, $z_0=3y_0$ is of infinite order.
\end{proof}
\begin{proof}[Proof of Theorem $\ref{th:1}$]
We note that
\[\sum_{d}z_d=3^kz_0,\]
where $d$ runs through all $d=p_1^{\epsilon_1}p_2^{\epsilon_2}\cdots p_k^{\epsilon_k}$ with $\epsilon_i=-1,0,1$. We will examine for which $d$ the terms $z_d$
will contribute on the left hand side.
\par If $k=1$, using Proposition $\ref{pro:4}$ and Proposition $\ref{pro:5}$, one can easily compute that $z_{1}=z_{{p_1^{* -1}}}=0$ and $z_{{p_1^*}}=3z_0$.
\par Assume $k\geq 2$. Let $d={p_1^{*}}^{\epsilon_1}p_2^{*\epsilon_2}\cdots p_k^{*\epsilon_k}$, where $\epsilon_i=-1,0,1$. Since $p_i\equiv 2\mod 3$, $1\leq i\leq k$, they are inert in $K$ and hence $E$ is supersingular at these primes. So the Fourier coefficients $a_{p_i}=0$ for the new form of level $\Gamma_0(36)$ and, by the Euler system property of CM points, see \cite{Nek}, if $\epsilon_i=0$ for some $i$, then $\Tr_{H_{6N}/H_{6N/p_i}}P_0=0$ and consequently $z_d=0$.
\par Since $y_0^{\sigma_{\omega_3}}=[\omega] y_0+t$, for some nonzero $t\in
E(\BQ)[3]$, we have
\[\left (\RTr_{H_0/K(\sqrt[3]{d})}y_0\right )^{\sigma_{\omega_3}}= [\omega] \RTr_{H_0/K(\sqrt[3]{d})}y_0\in E(K(\sqrt[3]{d}) .\]
By Proposition $\ref{pro:4}$,
\[\chi_d(\sigma_{\omega_3})=(\sqrt[3]{d})^{\sigma_{\omega_3}-1}=\omega^{\sum_{i=1}^k\epsilon_i}.\]
\par If $3\mid \sum_{i=1}^k\epsilon_i$, in which case the epsilon factor $\epsilon(E,\chi_d)=+1$, then $(\sqrt[3]{d})^{\sigma_{\omega_3}-1}=1$ and $\sigma_{\omega_3}=1$ on $K(\sqrt[3]{d})$. Then
\[\RTr_{H_0/K(\sqrt[3]{d})}y_0=\left (\RTr_{H_0/K(\sqrt[3]{d})}y_0\right )^{\sigma_{\omega_3}}=[\omega]\RTr_{H_0/K(\sqrt[3]{d})}y_0.\]
So $[1-\omega]\left (\RTr_{H_0/K(\sqrt[3]{d})}y_0\right )=0$ and therefore $\RTr_{H_0/K(\sqrt[3]{d})}z_0=3\RTr_{H_0/K(\sqrt[3]{d})}y_0=0$. Then we have
\[z_d=\sum_{\sigma\in \Gal(K(\sqrt[3]{d})/K)}[\chi_d^{-1}(\sigma)]\left (\RTr_{H_0/K(\sqrt[3]{d})}z_0\right )^\sigma=0.\]
\par If $3\nmid \sum_{i=1}^{k}\epsilon_i$, in which case the epsilon factor $\epsilon(E,\chi_d)=-1$, then $(\sqrt[3]{d})^{\sigma_{\omega_3}-1}\neq 1$ and $\Gal(K(\sqrt[3]{d})/K)=\langle\sigma_{\omega_3}\rangle$. We have
\[\left (\RTr_{H_0/K(\sqrt[3]{d})}z_0\right )^{\sigma_{\omega_3}}=[\omega]\RTr_{H_0/K(\sqrt[3]{d})}z_0.\] and then
\[z_d=\left \{
\begin{aligned}
3\RTr_{H_0/K(\sqrt[3]{d})}z_0 ,& \quad\sum_{i=1}^{k}\epsilon_i\equiv 1\mod 3;\\
0,\quad\quad\quad & \quad\sum_{i=1}^{k}\epsilon_i\equiv 2\mod
3.
\end{aligned}
   \right .\]
\par Thus we have an identity
\[3^kz_0=\sum_{\begin{subarray}{c}d={p_1^{*}}^{\epsilon_1}p_2^{*\epsilon_2}\cdots p_k^{*\epsilon_k}\\\forall i\,\,\epsilon_i\neq 0\\
\sum_1^k\epsilon_i\equiv 1\mod 3 \end{subarray}} z_d\] in $E(H_0)$. Here the condition $\sum_1^k\epsilon_i\equiv 1\mod 3$ exactly means that $\epsilon(1,E,\chi_d)=-1$ and $\epsilon(E^{(d)})=-1$. Since $z_0$ is of infinite order, at least one of such $z_d\in E(K(\sqrt[3]{d}))^{\chi_d}$ is of infinite order, and this produces a nonzero point in $E^{(d)}(K)_\BQ$. Then by the Gross-Zagier formula, $L(1,E^{(d^{-1})})\neq 0$ and hence $E^{(d^{-1})}(\BQ)$ is torsion and $2d^{-1}$ is not a cube sum.
\par Thus for arbitrary integer $k\geq 1$, and arbitrary primes $p_1,p_2,\cdots, p_k$ which are all congruent to $2,5 \mod 9$, there exists a cube-free integer $n$ with exactly prime factors $p_1,p_2,\cdots,p_k$ such that $2n$ is a cube sum (resp. not a cube sum). The theorem follows.
\end{proof}

\section{An Explicit Gross-Zagier Formula}

In this section, we generalize the construction of the Heegner point
$z_n$ in the cube sum problem to a very general framework. A height
formula for this Heegner point is then obtained as an application of
the variation of Gross-Zagier formula in \cite{CST} Theorem 1.6.

\s{\bf Local Theory.}
Let $F$ be a nonarchimedean
local field of characteristic zero. For $F$, denote by $\CO$ the ring
of integers in $F$, $\omega$ a uniformizer,
$\fp$ its maximal ideal and
$q$ the cardinality of its residue field. Let $B$ be a
quaternion algebra defined over $F$ with a quadratic $F$-subalgebra
$K$. Let $e$ be the ramification index of $K/F$ if $K$ is nonsplit.
Denote by $G$ the algebraic group $B^\times$ over $F$ and
also write $G$ for $G(F)$. Let $\pi$ be an irreducible admissible
representation of $G$ which is always assumed to be generic if
$G \cong \GL_2$. Denote by $\omega$ the central character of $\pi$.
Let $\chi$ be a character on $K^\times$ such that
$\chi|_{F^\times}\cdot\omega = 1$. Let $\CP(\pi,\chi)$ be the
functional space $\Hom_{K^\times}(\pi,\chi^{-1}))$. By a theorem
of Tunnell and Saito, the space $\CP(\pi,\chi)$ has dimension at most
one and equals one if and only if
\[\epsilon(\pi,\chi) = \chi\eta(-1)\epsilon(B).\]
Here, $\epsilon(\pi,\chi)$ is the Rankin-Selberg root number for
$\sigma \times \pi_\chi$ where $\sigma$ is
the Jacquet-Langlands lifting of $\pi$ to $\GL_2(F)$ and
$\pi_\chi$ is the representation on $\GL_2(F)$ constructed from
$\chi$ via Weil representation. The quadratic character $\eta$
corresponds to the extension $K/F$ via class field theory and
$\epsilon(B) = +1$ (resp. $=-1$) if $B$ is split (resp. division).

For our purpose, assume (1) the central character $\omega$ is
unramified (2) the pair $(\pi,\chi)$ is essentially unitary
in the sense that there exists some character $\mu = |\cdot|^s$ on
$F^\times$ with $s \in \BC$ such that both $\pi \otimes \mu$ and
$\chi\otimes \mu_K^{-1}$ are unitary. Here for a character $\mu$,
$\mu_K$ is the composition of $\mu$ with the norm map from $K$ to
$F$. Modify $(\pi,\chi)$ to $(\pi\otimes \mu, \chi\otimes \mu_K^{-1})$
for some $\mu$, we may from now on assume that $\pi$ and $\chi$ are
both unitary while $\pi$ has trivial central character.

Denote by $n$ the conductor of $\sigma$, that is,
the minimal non-negative integer $n$ such that the invariant subspace
of $\sigma^{U_0(n)}$ is nonzero where
\[U_0(n) = \left\{
	\begin{pmatrix}
		a & b \\
		c & d
	\end{pmatrix} \in \GL_2(\CO) \Big| c \in \fp^n
\right\}.\]
Let $c$ be the minimal non-negative integer such that
$\chi$ is trivial on $\CO_c^\times$
where $\CO_c = \CO + \varpi^c\CO_K$ and $\CO_K$
is the ring of integers for $K$.

\begin{lem}\label{epsilon}
	Let the pair $(\pi,\chi)$ be as above such that
	$\epsilon(\pi,\chi) = \chi\eta(-1)\epsilon(B)$. If
	$K$ is nonsplit and $c-n + e - 1 \geq 0$, then
	$B$ is split.
\end{lem}
\begin{proof}
  If $c \geq n$ or $n \geq 3$, then $B$ is split (see
  \cite{CST} Lemma 3.1 (5) and its proof). In the following, assume $K/F$ is ramified,
  $n = 2$ and $c = 1$. Assume $\pi$ is on $\GL_2(F)$ and  square-integrable.
  If $\pi = \mathrm{sp}(2) \otimes \mu$, then
  $\mu$ is a quadratic character on $F^\times$ with conductor $1$.
  By \cite{CST}, Lemma 3.1 (3), $B$ is nonsplit if and only if
  $\mu_K \cdot\chi = 1$. Since $\mu_K|_{F^\times} = 1$
  and the conductor of $\mu$ is $1$, $\mu_K$ is unramified. Hence,
  in this case, $B$ must be split. If $\pi$ is supercuspidal,
  then by \cite{Tunnell2}, Lemma 3.2, the two unramified characters
  in the dual of $K^\times/F^\times$ appear in $\pi_0$. Here,
  $\pi_0$ is the irreducible representation on the division algebra
  whose Jacquet-Langlands lifting is $\pi$. Since
  $\pi_0$ has dimension two, all the ramified characters should occur
  in $\pi$ and $B$ is split.
\end{proof}

Let $(\cdot,\cdot)$ be any non-degenerate invariant Hermitian pairing
on $\pi$. For any $f \in \pi$, denote by
\[\beta^0(f) = \int_{F^\times \bs K^\times} \frac{(\pi(t)f,f)}{(f,f)} \chi(t)dt\]
where $dt$ is any Haar measure on $F^\times \bs K^\times$. The
integral is absolutely convergence. The functional space
$\CP(\pi,\chi)$ is nontrivial if and only if $\beta^0$ is nontrivial.
If $\CP(\pi,\chi)$ is nonzero, a nonzero vector $f$ of $\pi$ is
called a {\em test vector} for $\CP(\pi,\chi)$ if $\ell(f) \not=0$
for some (thus any) nonzero $\ell \in \CP(\pi,\chi)$, or
equivalently, $\beta^0(f)$ is non-vanishing.

We shall end our discussion of local theory by the following lemma
which will be used in the proof of main result Theorem \ref{gz-sharp}
\begin{lem}\label{lemma-toric}
   Assume $K/F$ is ramified, $cn \not=0$ with $c+1 = n$. Let $B$
   be split. Let $R$ be an order in $B$ with discriminant $n$ such that
   $R  = R' \cap R''$. Here $R'$ and $R''$ are two maximal orders in $B$
   such that $R' \cap K = \CO_c$ and $R'' \cap K = \CO_K$. Such
   order exists and unique up to $K^\times$-conjugacy. The invariant
   subspace $\pi^{R^\times}$ is of dimension $1$. For any
   nonzero $f \in \pi^{R^\times}$,
   \[\beta^0(f) = 2^{-1}q^{-c}\Vol(K^\times/F^\times)L(1,1_F).\]
\end{lem}
\begin{proof}
 The existence and uniqueness of the order $R$ follows from \cite{CST} Lemma 3.2.
 The order $R$ is an Eichler order. In particular, if we fix an
 isomorphism $B \cong M_2(F)$, then $R^\times$ is $\GL_2(F)$-conjugate
 to $U_0(n)$. By new form theory, the invariant subspace $\pi^{R^\times}$ is
 of dimension one. Next, we compute the local toric integral $\beta^0(f)$ for
 any $f$. It is similar to that in \cite{CST} \S 3.4.
     Let $\tau$ be a uniformizer of $K$
     such that $\CO_K = \CO[\tau]$ and embed $K$ into $M_2(F)$ by
     \[a+b\tau \mapsto \gamma_c^{-1}
       \begin{pmatrix}
	 a+b\mathrm{tr}\tau & b\mathrm{N}\tau \\
	 -b & a
       \end{pmatrix}\gamma_c, \quad \gamma_c =
       \begin{pmatrix}
	 \varpi^c\mathrm{N}\tau & \\
	   & 1
       \end{pmatrix},\]
       where $\omega$ is a uniformizer of $F$ and
       $\ord_v(c_v)$ is denoted by $c$.
       In particular, under this embedding, $R_0(n) =
       \begin{pmatrix}
	 \CO & \CO \\
	 \fp^{n} & \CO
       \end{pmatrix}$  is an admissible order for $(\pi,\chi)$. Since
       the toric integral is invariant under $K^\times$-conjugacy of
       $f$, we may take $R = R_0(n)$ and $f=f_0$ the normalized
       new form of $\pi$. Let $\Psi(g)$ denote the matrix coefficient
\[\Psi(g) := \frac{(\pi(g)f_0,f_0)}{(f_0,f_0)}, \qquad  g\in \GL_2(F).\]
Then $\Psi$ is a function on $ZU_0(n) \bs \GL_2(F) / U_0(n)$ where
$Z$ is the center of $\GL_2(F)$ and
\[\beta(f)
    =\frac{\Vol(K^\times/F^\times)}{\# K^\times/F^\times\CO_c^\times}
  \sum_{t \in K^\times/F^\times\CO_c^\times} \Psi(t)\chi(t).\]
Denote by
    \[S_i = \left\{ 1+b\tau, b \in \CO/\fp^c, v(b) = i \right\}, \quad 0 \leq i \leq c-1\]
    and
    \[S' =\left\{ a\varpi+\tau, a \in \CO/\fp^c \right\}.\]
Then a complete representatives of $K^\times/F^\times\CO_c^\times$ can be taken as
    \[\left\{ 1 \right\}\sqcup\left(\sqcup_{i} S_i\right)\sqcup S'.\]
    Let $\mathrm{pr}: K^\times \ra ZU_0(n) \bs \GL_2(F) / U_0(n)$ be the natural mapping. Then
it is constant on $S_i$ and $S'$. Precisely,
\[\mathrm{pr}(S_i) = \left[
    \begin{pmatrix}
      1 & \varpi^{i-c} \\
        & 1
    \end{pmatrix}
      \right], \quad
      \mathrm{pr}(S') = \left[
    \begin{pmatrix}
      & \varpi^{-c} \\
      -\varpi^{c+1} &
    \end{pmatrix}
      \right].\]
Follow from this
    \[\sum_{t \in K^\times/F^\times\CO_c^\times} \Psi(t)\chi(t) =
    1+\sum_{i=0}^{c-1} \Psi_i \sum_{t\in S_i} \chi(t) + \Psi'\sum_{t \in S'} \chi(t),\]
where $\Psi_i$ (resp. $\Psi'$) are the valuations of $\Psi(t)$ on $S_i$ (resp. $S'$). As
      \[\sum_{t \in S_i} \chi(t) =
    \begin{cases}
      0, \quad &\text{if $c>1$ and } 0 \leq i \leq c-2, \\
      -1,\quad &\text{if } i=c-1,
      \end{cases} \quad
      \text{and} \quad
      \sum_{t \in S'} \chi(t) = 0,\]
    we obtain,
    \[ \sum_{t \in K^\times/F^\times\CO_c^\times} \Psi(t)\chi(t) = 1 - \Psi_{c-1}.\]
Finally, one needs to compute this matrix coefficient.  However,
as it actually evaluates at some upper-triangler matrix, the computation
is easy by the explicit description of normalized new forms in
Kirillov model. In fact, $\Psi_{c-1} = -q^{-1}L(1,1_F)$ where
$q$ is the cardinality of the residue field of $F$ (See also \cite{CST} \S 3.4).
Thus
\[\beta^0(f) = \frac{\Vol(K^\times/F^\times)}
{\# K^\times/F^\times\CO_{c}^\times}L(1,1_F) = 2^{-1}q^{-c}\Vol(K^\times/F^\times)L(1,1_F).\]
\end{proof}

\s{\bf Global Theory.}
Let $F$ be a totally real field with $\CO$ its ring of integers.
Take $d := [F: \BQ]$.
Denote by $\BA = F_\BA$ be the adele ring of $F$ and $\BA_f$ its
finite part. For a $\BZ$-module $M$, take $\wh{M} = M \otimes_\BZ \wh{\BZ}$
with $\wh{\BZ} = \prod_{p<\infty} \BZ_p$.
Let $\BB$ be a totally
definite incoherent quaternion algebra over $\BA$. For any open compact subgroup
$U$ of $\BB_f^\times = (\BB \otimes_\BA \BA_f)^\times$,
denote by $X_U$ the Shimura curve associated to $\BB^\times$ of level $U$.
Denote by $\xi_U \in \Pic(X_U)_\BQ$ the normalized Hodge class on $X_U$.
Let $X$ be the projective limit of the system $(X_U)_U$.
Let $A$ be a simple abelian variety defined over $F$ parametrized by $X$
in the sense that there is a non-constant morphism $X_U \ra A$ over $F$
for some $U$. Then $M:=\End^0(A)$ is a field.
Denote by $\pi = \pi_A = \varinjlim_U \Hom^0_{\xi_U} (X_U,A)$.
Here, $\Hom_{\xi_U}^0(X_U,A)$ denotes the morphisms in
$\Hom(X_U,A)\otimes_\BZ \BQ$ using $\xi_U$ as a base point: if
$\xi_U$ is represented by a divisor $\sum_i a_i x_i$ on $X_{U,\bar{F}}$,
then $f \in \Hom_F(X_U,A)\otimes_\BZ \BQ$ is in $\pi_A$ if and only if
$\sum a_i f_i(x_i) = 0$ in $A(\bar{F})_\BQ = A(\bar{F}) \otimes_\BZ \BQ$.
For each open compact subgroup $U$ of $\BB_f^\times$, let $J_U$ denote
the Jacobian of $X_U$. Then $\pi = \varinjlim_U \Hom^0(J_U,A)$ where
$\Hom^0(J_U,A) = \Hom_F(J_U,A)\otimes_\BZ \BQ$.
The Hecke action of $\BB^\times$ on $X$ induces a
natural $\BB^\times$-module structure on $\pi$ so that
$\End_{\BB^\times}(\pi) = M$ and has a decomposition
$\pi = \otimes_M \pi_{v}$ where $\pi_v$ are absolutely irreducible
representations of $\BB_v^\times$ over $M$.
Let $N$ be the conductor of the Jacquet-Langlands lifting of $\pi$.
Denote by $\omega_A$ its central character.
Let $A^\vee$ be the dual abelian variety of $A$. There is a perfect
$\BB^\times$-invariant pairing
\[\pi_A \otimes \pi_{A^\vee} \lra M\]
given by
\[(f_1,f_2) = \Vol(X_U)^{-1}(f_{1,U} \circ f_{2,U}^\vee),
\quad f_{1,U} \in \Hom(J_U,A), f_{2,U} \in \Hom(J_U,A^\vee)\]
where $f_{2,U}^\vee: A \ra J_U$ is the dual of $f_{2,U}$ composed with
the canonical isomorphism $J_U^\vee \cong J_U$. Here,
$\Vol(X_U)$ is defined by a fixed invariant measure on the upper half
plane.
It follows that $\pi_{A^\vee}$ is dual to $\pi_A$ as representations
of $\BB^\times$ over $M$.
For any fixed open compact subgroup $U$ of $\BB_f^\times$,
define the $U$-pairing on $\pi_A \times \pi_{A^\vee}$ by
\[(f_1,f_2)_U = \Vol(X_U)(f_1,f_2),
\quad f_1 \in \pi_A, f_2 \in \pi_{A^\vee}\]
which is independent of the choice of measure defining $\Vol(X_U)$.

Let $K$ be a totally imaginary quadratic extension over $F$ with relative discriminant $D = D_{K/F} \subset \CO$
and absolute discriminant $D_K$.
Denote by $\CO_K$ its ring of integers. Let $\eta$ be the quadratic character corresponding to $K/F$.
Take $\chi: K^\times \bs K_\BA^\times \ra L^\times$ a character of finite order
where $L$ is a finite extension on $M$.
Let $c \subset \CO$ be the ideal maximal such that $\chi$ is
trivial on $\prod_{v \not| c} \CO_{K,v}^\times \prod_{v|c} (1+c\CO_{K,v})$.
Assume that $\omega_A\cdot\chi|_{\BA_f^\times} = 1$
and for any place $v$,
\[\epsilon(\pi_v,\chi_v) = \chi_v\eta_v(-1)\epsilon(\BB_v),\]
where $\epsilon(\BB_v) = 1$ if $\BB_v$ is split or $=-1$ otherwise
and $\epsilon(\pi_v,\chi_v)$ is the local root number of
the Rankin-Selberg $L$-series $L(s,\pi,\chi)$.
These assumptions imply that the Rankin-Selberg root number
$\epsilon(\pi,\chi) = -1$ and there exists an embedding
$K_\BA^\times \ra \BB^\times$.
Let $X^{K^\times}$ be the $F$-subscheme of $X$ of fixed points of $X$ under
$K^\times$. The theory of complex multiplication asserts that every
point in $X^{K^\times}(\bar{F})$ is defined over $K^{\mathrm{ab}}$ and
that the Galois action is given by the Hecke action under the reciprocity
law. Let $P \in X^{K^\times}(K^{\mathrm{ab}})$. For any $f \in \pi$, consider the Heegner cycle
\[\int_{K^\times \bs \wh{K}^\times/\wh{F}^\times} f(P)^{\sigma_t}\chi(t)dt \in A(K^{\mathrm{ab}})_\BQ \otimes_M L\]
where $\sigma: K^\times \bs \wh{K}^\times/\wh{F}^\times \ra \Gal(K^{\mathrm{ab}}/K)$ is the Artin map. The
main result in the book of Yuan-Zhang-Zhang \cite{YZZ} implies that there is some
$f \in \pi$ such that the corresponding Heegner cycle for $f$ is nonvanishing  if and only if  $L'(1,A,\chi)$ is nonvanishing.

In the following, we shall explicit construct a one-dimensional $L$-subspace $V(\pi,\chi)$
of $\pi\otimes_M L$ such that for any nonzero $f \in V(\pi,\chi)$, the Heegner cycle for $f$ is nonzero
if and only if $L'(1,A,\chi)$ is nonzero. A similar one-dimensional $L$-subspace is
considered in \cite{CST} which is also denoted by $V(\pi,\chi)$ there. Any vectors in these two
$L$-lines are of pure tensor form. In fact, any two vectors from these two $L$-lines respectively
are parallel at any places
outside
\[\BS := \left\{ v|(c,N) \Big|
   \text{$K_v/F_v$ is ramified and $\ord_v(c) + 1= \ord_v(N)$} \right\}.\]

From now on, we require that the central character $\omega$ is unramified at any place in
$\BS$. Define the following sets of places $v$ of $F$ dividing $N$:
\[\Sigma_1 := \left\{ v | N \text{ nonsplit in $K$}: \ord_v(c) = 0 \text{ or }
   \ord_v(c) - \ord_v(N) + e_v -1 <0 \right\},\]
where  $e_v$ is the ramification index of $K$ at $v$.
Let $c_1=\prod_{\fp|c, \fp\notin \Sigma_1} \fp^{\ord_\fp c}$
be the $\Sigma_1$-off part of $c$,  $N_1$ the $\Sigma_1$-off part of $N$,
and $N_2=N/N_1$. Denote by $\CO_{c_1} = \CO + c_1\CO_K$.

Let $v$ be a place of $F$ and $\varpi_v$ a uniformizer of $F_v$.
Then there exists an $\CO_v$-order $R_v$ of $\BB_v$ with discriminant $N\CO_v$
such that $R_v \cap K_v=\CO_{c_1, v}$.
Such order exists and unique up to  $K_v^\times$-conjugate  when $v\nmid (c_1, N)$.
At any place $v|(c_1,N)$, by Lemma \ref{epsilon}, $\BB_v$ is split and $R_v$ will
be taken as an Eichler order. In particular, fixing an isomorphism
between $\BB_v$ and $M_2(F_v)$, then $R_v$ is $\GL_2(F_v)$-conjugate to
the order
$\begin{pmatrix}
   \CO_v & \CO_v \\
   N\CO_v & \CO_v
\end{pmatrix}$. However, the $K_v^\times$-conjugacy class of $R_v$ needs to be
determined here. Such an order $R_v$ is called {\em admissible} for $(\pi_v, \chi_v)$ if it also satisfies
the following  conditions (1) and (2).
\begin{enumerate}
  \item If $v|(c_1, N)$,
    then $R_v$ is the intersection of two maximal orders $R_v', R_v''$ of $\BB_v$ such that
    $R_v'\cap {K_v}=\CO_{c, v}$ and  $$ R_v''\cap K_v=\begin{cases}\CO_{c/N, v},
     \  &\text{ if $\ord_v(c/N)\geq 0$},\\
       \CO_{K, v}, &\text{otherwise.}\end{cases}$$
   \end{enumerate}
Note that for $v|(c_1, N)$, By \cite{CST} Lemma 3.2, there is a unique order,
up to $K_v^\times$-conjugate,
satisfying the condition (1) unless $K_v$ is split and $0<\ord_v(c_1)<\ord_v(N)$.
In the case  $K_v$ is split and $0<\ord_v(c_1)<\ord_v(N)$,
there are exactly two $K_v^\times$-conjugacy classes of orders satisfying the condition (1),
which are conjugate to each other by a normalizer of $K_v^\times$ in $\BB_v^\times$.
Fix an $F_v$-algebra isomorphism $K_v\cong F_v^2$ and identify $\BB_v$ with $\End_{F_v}(K_v)$.
Then the two classes contain respectively orders $R_{i, v}=R'_{i, v}\cap R''_{i, v}, i=1, 2$ as in (1)
such that $R'_{i, v}=\End_{\CO_v} (\CO_{c,v}), i=1, 2,$ and
$R''_{1, v}=\End_{\CO_v}((\varpi_v^{n-c}, 1)\CO_{K_v})$ and
$R''_{2, v}=\End_{\CO_v}((1, \varpi_v^{n-c})\CO_{K_v})$.
\begin{enumerate}
  \item[(2)] If $K_v$ is split and $0<\ord_v(c_1)<\ord_v(N)$,
    then  $R_v$ is  $K_v^\times$-conjugate to some $R_{i, v}$
    such that $\chi_i$ has conductor $\ord_v(c)$,
    where $\chi_i, i=1, 2$ is defined by $\chi_1(a)=\chi_v(a, 1)$ and $\chi_2(b)=\chi_v(1, b)$.
  \end{enumerate}
\begin{defn}
An $\wh{\CO}$-order $\CR$ of $\BB_f$ is called admissible for $(\pi, \chi)$
if for every finite place $v$ of $F$, $\CR_v:=\CR\otimes_{\wh{\CO}}\CO_v$ is
admissible for $(\pi_v, \chi_v)$.
Note that an admissible order $\CR$ for $(\pi, \chi)$ is of discriminant $N\wh{\CO}$
such that $\CR\cap \wh{K}=\wh{\CO}_{c_1}$.
\end{defn}

Let $\CR$ be an $\wh{\CO}$-order of $\BB_f$ which is admissible for $(\pi,\chi)$.
Let $U=\CR^\times$ and $U^{(N_2)}:=\CR^\times \cap \BB_f^{\times (N_2)}$.
Note that for any finite place $v|N_1$, $\BB_v$ must be split.
Let $Z\cong \BA_f^\times$ denote the center of $\BB_f^\times$.
The group $U^{(N_2)}$ has a decomposition $U^{(N_2)}=U' \cdot(Z\cap U^{(N_2)})$
where $U'=\prod_{v\nmid N_2\infty} U'_v$
such that for any finite place $v\nmid N_2$,
$U'_v=U_v$ if $v\nmid N$ and $U'_v\cong {U_1(N)}_v$ otherwise.
View $\omega$ as a character on $Z$ and we may define a character on $U^{(N_2)}$
by $\omega$ on $Z\cap U^{(N_2)}$ and trivial on $U'$, which we also denoted by $\omega$.

\begin{defn} Let $V(\pi, \chi)$ denote the space of forms $f\in \pi_A\otimes_M L$,
  which are $\omega$-eigen under $U^{(N_2)}$,
  and $\chi_v^{-1}$-eigen under $K_v^\times$ for all places $v\in \Sigma_1$.
  The space $V(\pi, \chi)$ is actually a one dimensional $L$-space.
\end{defn}

For any $f \in V(\pi,\chi)$, define
the Heegner cycle associated to $(\pi,\chi)$ to be
\[P_\chi(f) := \sum_{t \in \wh{K}^\times/K^\times\wh{F}^\times\wh{\CO}_{c_1}^\times} f(P)^{\sigma_t}\chi(t)
 \in A(K^{\mathrm{ab}})_\BQ \otimes_M L.\]

The Neron-Tate height pairing over $K$ gives a $\BQ$-linear map $\langle , \rangle_K:
      A(\bar{K})_\BQ \otimes_M A^\vee(\bar{K})_\BQ \ra \BR$. Let $\langle , \rangle_{K,M}:
      A(\bar{K})_\BQ \otimes_M A^\vee(\bar{K})_\BQ \ra M\otimes_\BQ \BR$ be the unique $M$-bilinear
      pairing such that $\langle , \rangle_K = \tr_{M\otimes \BR/\BR} \langle , \rangle_{K,M}$. This
      induces a $L$-linear Neron-Tate pairing over $K$:
      \[\left\langle \cdot,\cdot \right\rangle_{K,L}:
	(A(\bar{K})_\BQ \otimes_M L) \otimes (A^\vee(\bar{K})_\BQ \otimes_M L)
      \lra L\otimes_\BQ \BR.\]
Let $\phi$ be the Hilbert new form in the Jacquet-Langlands lifting of $\pi$ defined as following.
For $v|\infty$, the subgroup $\SO_2(\BR)$ of $\GL_2(F_v)$ acts on $\phi$ via the character $k_\theta =
    \begin{pmatrix}
      \cos\theta & \sin\theta \\
      -\sin\theta & \cos\theta
    \end{pmatrix} \mapsto e^{4\pi i\theta}$. It is of level $U_1(N)$ with
    $$ U_1(N)=\left\{ \matrixx{a}{b}{c}{d}\in \GL_2(\wh{\CO})\Bigg|
    \begin{aligned} &c \in N\wh{\CO} \\ &d\equiv 1\mod N\wh{\CO} \end{aligned} \right\}.$$
Finally,
    \[L(s,\pi) = 2^d |\delta|_\BA^{s-1/2}Z(s,\phi), \quad Z(s,\phi) = \int_{F^\times\bs \BA^\times}
      \phi\left[
	\begin{pmatrix}
	  a & \\
	    & 1
	\end{pmatrix}
      \right]|a|_\BA^{s-1/2} d^\times a\]
where the Haar measure $d^\times a$ is chosen such that $\Res_{s=1} \int_{|a| \leq 1} |a|^{s-1}d^\times a  = \Res_{s=1} L(s,1_F)$.
The Petersson norm $(\phi,\phi)_{U_0(N)}$ is the integration of $\phi\bar{\phi}$ with respect to the measure $dxdy/y^2$
on the upper half plane.

\begin{thm}\label{gz-sharp} Let $F$ be a totally real field of
	degree $d$. Let $A$ be an abelian variety over $F$ parametrized by
	a Shimura curve $X$ over $F$ associated to an incoherent totally
	definite quaternion algebra $\BB$ over $\BA$. Denote by
	$\phi$ the Hilbert holomorphic
	newform of parallel weight $2$ on $\GL_2(\BA)$ associated to $A$.
	Let $K$ be a totally imaginary quadratic extension over $F$ with
	relative disriminant $D$ and discriminant $D_K$. Let
	$\chi: K_\BA^\times/K^\times \ra L^\times$ be a finite Hecke
	character of conductor $c$ over some finite extension $L$
	of $M:= \End^0(A)$. Assume that
	\begin{enumerate}
		\item $\omega_A\cdot \chi|_{\BA^\times} = 1$,
			where $\omega_A$ is the central character of $\pi_A$;
	        \item for any place $v$ of $F$,
			$\epsilon(\pi_{A,v},\chi_v) = \chi_v\eta_v(-1)
			\epsilon(\BB_v)$;
	        \item $\omega$ is unramified at any place in $\BS$.
	\end{enumerate}
  For any non-zero forms $f_1 \in V(\pi_A,\chi)$ and
  $f_2 \in V(\pi_{A^\vee},\chi^{-1})$,
  we have an equality in $L\otimes_\BQ \BC$:
  \[L^{'(\Sigma)}(1,A,\chi) = 2^{-\#(\Sigma_D)}\frac{(8\pi^2)^d(\phi,\phi)_{U_0(N)}}
    {(u_1)^2\sqrt{|D_K|||c_1||^2}}\cdot
  \frac{\langle P_\chi(f_1),P_{\chi^{-1}}(f_2)\rangle_{K,L}}{(f_1,f_2)_{\CR^\times}}.\]
    In the above,
  \[\Sigma = \left\{ v|(N,Dc) \big| \text{if $v||N$ then  $\ord_v(c/N) \geq 0$} \right\}\]
  and
  \[\Sigma_D = \left\{ v|(N,D) \big| v \not| c \text{ or } 0<\ord_v(c) < \ord_v(N) -1 \right\}.\]
  The term $c_1$ is the $\Sigma_1$-off part of $c$ and
  $u_1 = \#\ker(\Pic(\CO) \ra \Pic(\CO_{c_1}))\cdot [\CO_{c_1}^\times:\CO^\times]$ while
  $||c_1||$ is the norm of $c_1$.
\end{thm}
\begin{proof}
  %By the variation of Gross-Zagier formula in \cite{CST} Theorem 1.6, the proof reduces to compare
  %the local toric integrals between $V(\pi,\chi)$ and $V(\pi,\chi)^\sharp$ at any place in
  %$\Sigma^\sharp$. Similar computations as in \cite{CST} obtains that
  %for any $f \in V(\pi,\chi)^\sharp$ and $v \in \Sigma^\sharp$,
  %the local toric integral $\beta(f_v)|D_v\delta_v|^{-1/2}$ equals
  %\[\frac{L(1,1_v)}{L(2,1_v)}L(1,\eta_v)^2|\varpi_v^{c_v}|
  %\frac{L(1,\pi_v,\ad)}{L(1/2,\pi_v,\chi_v)}.\]
  %Note that this value is the same as that for $V(\pi,\chi)$
  %under the case $c_v \geq n_v > 0$ in \cite{CST} Proposition 3.12.

  The proof is an application of the variation of Gross-Zagier formula (\cite{CST} Theorem 1.6).
  In other words, the above height formula is derived from the explicit Gross-Zagier formula
  in \cite{CST} Theorem 1.5 plus some local computations. In the following proof, we shall
  use notations in \cite{CST} Theorem 1.5 and add the sign $\sharp$ for any corresponding
  terms in Theorem \ref{gz-sharp}. Moreover, we shall also use measures defined in \cite{CST}.

  By \cite{CST} Theorem 1.6, for any $f_1^\sharp \in V^\sharp(\pi_A,\chi)$ and
  $f_2^\sharp \in V^\sharp(\pi_{A^\vee},\chi^{-1})$,
  \[\prod_{v \in \BS} \frac{\beta^0(f_{1,v}^\sharp,f_{2,v}^\sharp)}{\beta^0(f_{1,v},f_{2,v})}
    L^{'(\Sigma)}(1,A,\chi) = 2^{-\# \Sigma_D}\cdot \frac{ (8\pi^2)^d (\phi,\phi)_{U_0(N)}}
    {u_1^2\sqrt{|D_K|||c_1||^2}}\cdot \frac{ \langle P_\chi^{c_1}(f_1^\sharp),
  P_{\chi^{-1}}^{c_1}(f_2^\sharp)\rangle_{K,K}}{(f_1^\sharp,f_2^\sharp)_{\CR^\times}}.\]

  The definition of $\Sigma$ is the same as the one in Theorem \ref{gz-sharp}.
    On the other hand,
    \[\Sigma_D = \{ v| (N,D) | \ord_v(c) < \ord_v(N)\} = \Sigma_D^\sharp \sqcup \BS.\]

  The term $c_1$ is the $\Sigma_1$-off part of $c$. Note that in
    \cite{CST}, $\Sigma_1$ is defined as the set of places $v|N$ nonsplit
    in $K$ with $\ord_v(c) \leq \ord_v(N)$. In particular, if
    we denote by $\Sigma_1^\sharp$ the corresponding set defined before Lemma \ref{epsilon},
    then $\Sigma_1 = \Sigma_1^\sharp \sqcup \BS$. Thus, $c_1c_\BS = c_1^\sharp$ and
    $||c_1|| = |c|_\BS||c_1^\sharp||$. The Heegner point
    \[P^{c_1}_\chi(f_1^\sharp) = \frac{\# \Pic_{K/F}(\CO_{c_1})}{\# \Pic_{K/F}(\CO_{c_1^\sharp})}
    \sum_{t \in \Pic_{K/F}(\CO_{c_1^\sharp})} f_1^\sharp(P)^{\sigma_t}\chi(t)
     = \frac{\# \Pic_{K/F}(\CO_{c_1})}{\# \Pic_{K/F}(\CO_{c_1^\sharp})}\cdot P_\chi(f_1^\sharp)\]
    and the definition of $P^{c_1}_{\chi^{-1}}(f_2^\sharp)$ is similarly.
    Denote by $\kappa_{c_1} = \ker(\Pic(\CO) \ra \Pic(\CO_{c_1}))$ and similarly for
    $\kappa_{c_1^\sharp}$, then by \cite{CST} Lemma 2.3,
    \[
      \quad\quad\quad\frac{\# \Pic_{K/F}(\CO_{c_1})}{\# \Pic_{K/F}(\CO_{c_1^\sharp})}
      = \frac{[\wh{\CO}_K^\times: \wh{\CO}_{c_1}^\times]}{[\wh{\CO}_K^\times:\wh{\CO}_{c_1^\sharp}^\times]}
      \frac{[\CO_K^\times:\CO_{c_1}^\times]^{-1}}{[\CO_K^\times:\CO_{c_1^\sharp}^\times]^{-1}}
    \frac{\#\kappa_{c_1}}{\#\kappa_{c_1^\sharp}}
      = \prod_{v \in \BS} \frac{1}{[\CO_{K,v}^\times:\CO_{c,v}^\times]}
      [\CO_{c_1}^\times:\CO_{c_1^\sharp}^\times]
      \frac{\#\kappa_{c_1}}{\#\kappa_{c_1^\sharp}} = |c|_\BS\frac{u_1}{u_1^\sharp}.
    \]

    The admissible order $\CR = \prod_v \CR_v$ is different with the $\sharp$-admissible
     order $\CR = \prod_v \CR_v^\sharp$ exactly at $v \in \BS$. For $v \in \BS$, $\CR_v$ is
     an order in $\BB_v$ with discriminant $\ord_v(N)$ such that $\CR_v \cap K_v = \CO_{K,v}$.
     By Lemma 2.2 and Lemma 3.5 in \cite{CST},
     \[\frac{\Vol(X_{\CR^\times})}{\Vol(X_{\CR^{\sharp \times}})} =
       \frac{\Vol(\CR^{\sharp \times})}{\Vol(\CR^\times)} = \prod_{v \in \BS}
     \frac{\Vol(\CR^{\sharp \times}_v)}{\Vol(\CR^\times_v)} = L(1,1_F)_\BS^{-1}.\]
     Thus
     \[(f_1^\sharp,f_2^\sharp)_{\CR^{\times}} = \frac{\Vol(X_{\CR^\times})}{\Vol(X_{\CR^{\sharp \times}})}
       (f_1^\sharp,f_2^\sharp)_{\CR^{\sharp\times}} = L(1,1_F)_\BS^{-1}
     (f_1^\sharp,f_2^\sharp)_{\CR^{\sharp\times}}.\]

  For any place $v \in \BS$, fix a $\BB_v^\times$-invariant pairing $\langle \cdot,\cdot \rangle_v
  $ on $\pi_{A,v} \times \pi_{A^\vee,v}$.  For any $f'_{1,v} \in \pi_{A,v}$ and $f'_{2,v} \in \pi_{A^\vee,v}$ with $
  \langle f'_{1,v}, f'_{2,v} \rangle_v \not=0$, let
  \[\beta^0(f'_{1,v},f'_{2,v}) = \int_{F_v^\times \bs K_v^\times}
  \frac{\langle \pi_{A,v}(t_v)f'_{1,v},f'_{2,v} \rangle_v}{\langle f'_{1,v},f'_{2,v} \rangle_v}\chi_v(t_v)dt_v.\]
  The test vector $f_{1,v}$ (resp. $f_{2,v}$) is the $v$-component of a nonzero vector in
  $V(\pi_A,\chi)$ (resp. $V(\pi_{A^\vee},\chi^{-1})$). As $v \in \BS$, $f_{1,v}$ (resp.
  $f_{2,v}$) is $\chi_v^{-1}$-eigen (resp. $\chi_v$-eigen) under $K_v^\times$.
  \begin{lem}
    For any $v \in \BS$,
    \[\frac{\beta^0(f_{1,v},f_{2,v})}{\beta^0(f_{1,v}^\sharp,f_{2,v}^\sharp)} = 2L(1,1_v)^{-1}|c_v|_v^{-1}.\]
  \end{lem}
  \begin{proof}
     As this is a local computation, we shall drop the subscript '$v$' in the
     following.
     The toric integral $\beta^0$ is invariant by modifying
     $(\pi,\chi)$ to $(\pi\otimes\mu,\chi\otimes\mu_K^{-1})$ for any
     character $\mu$ of $F^\times$. Since the central character $\omega$
     of $\pi$ is assumed to be unramified, we may assume $\omega = 1$.
     Identify the contragredient representation
     $(\pi^\vee,\chi^{-1})$ with the complex-conjugation
     $(\bar{\pi},\bar{\chi})$ and let $(\cdot,\cdot)$ be the
     Hermitian pairing on $\pi$ defined by $(f_1,f_2) = \langle f_1,\ov{f_2}
     \rangle$. Denote by $\beta^0(f) = \beta^0(f,\bar{f})$. Then
     the ratio we need to compute becomes
     $\frac{\beta^0(f)}{\beta^0(f^\sharp)}$ with $f \in V(\pi,\chi)$ and
     $f^\sharp \in V^\sharp(\pi,\chi)$.  Since $f \in V(\pi,\chi)$ is
     $\chi^{-1}$-invariant, $\beta^0(f) = \Vol(K^\times/F^\times)$. On the other hand,
     by Lemma \ref{lemma-toric}, $\beta^0(f^\sharp) = 2^{-1}q^{-c}\Vol(K^\times/F^\times)L(1,1_F)$.
     The result is then obtained.
  \end{proof}

 Sum up,
 \[L^{'(\Sigma)}(1,A,\chi) = 2^{-\#(\Sigma_D \setminus \BS)}\cdot\frac{(8\pi^2)^d(\phi,\phi)_{U_0(N)}}
   {(u_1^\sharp)^2\sqrt{|D_K|||c_1^\sharp||^2}}
   \frac{\langle P_\chi(f_1^\sharp),P_{\chi^{-1}}(f_2^\sharp)\rangle_{K,K}}
 {(f_1^\sharp,f_2^\sharp)_{\CR^{\sharp\times}}}.\]
 Using notations in Theorem \ref{gz-sharp}, this is just the formula we need to prove.
\end{proof}

\s{\bf Proof of Theorem $\ref{th:2}$.}
Let $n={p_1^{*}}^{\epsilon_1}p_2^{*\epsilon_2}\cdots p_k^{*\epsilon_k}$
with $\epsilon_i=\pm 1$ and $\sum_1^k\epsilon_i\equiv 1\mod 3$,
which is equivalent to $\epsilon(E,\chi_n)=-1$ and $\epsilon(E^{(n)})=-1$.
Denote by $N = p_1\cdots p_k$.
Let $f$ be the identity morphism from $X_0(36)$ to $E$ taking the cusp $[\infty]$ to the identity point of $E$.
Let $K = \BQ(\sqrt{-3})$ be the complex multiplication field of $E$.
Let $\chi_n:\Gal(H_{6N}/K)\rightarrow \BC^\times$ be the cubic character that $\chi_n(\sigma)=\sqrt[3]{n}^{\sigma-1}$
for $\sigma\in \Gal(H_{6N}/K)$. The Heegner point for $\chi_n$ is defined as
\[z_n=\sum_{\sigma\in \Gal(H_{6N}/K)}[\chi_n^{-1}(\sigma)]f(P_0^{\sigma}).\]

\begin{lem}\label{test-vector}
	Let $B$ be the incoherent quaternion algebra over $\BQ$ determined by the pair $(E,\chi_n)$. Then
	$B$ is split.  Let $\pi_E$ be the cuspidal automorphic representation associated to $E$, then
	$f \in V(\pi_E,\chi_n)$.
\end{lem}
\begin{proof}
   Note that the conductor of $\pi_E$ is $36$, the conductor $c(\chi_n)$ of $\chi_n$ is $3N$. If $p \nmid 36$, then
  $B_p$ is split by Lemma 3.1 (1) in \cite{CST}. As $2$ is inert in $K$, $B_2$ is split by Lemma 3.1 (4) in \cite{CST}.
  Finally, $B_3$ is also split by Lemma \ref{epsilon}. Hence $B$ is split. We next check that $f \in V(\pi_E,\chi_n)$.
  Note that by our choice of embedding of $K$ to $M_2(\BQ)$, $R_0(36) \cap K = \CO_{6N}$ with
  $R_0(36) =
  \begin{pmatrix}
	  \BZ & \BZ \\
	  36\BZ & \BZ
  \end{pmatrix}$. It suffices to check that $f$ is $K_2^\times$-invariant and $R_0(36)_3$ is admissible for $(\pi_{E,3},\chi_{n,3})$.
  The fact that $f$ is $K_2^\times$-invariant follows from that $K_2^\times/\BQ_2^\times(1+2\CO_{K,2})$ is
  generated by $\omega_2$ and by Theorem \ref{th:4}, for any point $P$ in $X_U$, $P^{\omega_2} = P + \tau(2)$.
  Finally, $R_0(36)_3 = R' \cap R''$ where $R' = M_2(\BZ_3)$ and
  $R'' =
  \begin{pmatrix}
	  \BZ_3 & 9^{-1}\BZ_3 \\
	  9\BZ_3 & \BZ_3
  \end{pmatrix}$. It is then easy to check that $R' \cap K_3 = \BZ_3 + 3\CO_{K,3}$ while
  $R'' \cap K_3 = \CO_{K,3}$. Thus, $R_0(36)_3$ is admissible for $(\pi_{E,3},\chi_{n,3})$.
\end{proof}

Let $\Omega^{(n)}$ denote the minimal real period of $E^{(n)}$ and $\Omega$ the one of $E$.
Then $\Omega^{(n)}\Omega^{(n^{-1})}=\Omega^2/N$ with $N=p_1\cdots p_k$.
Let $\phi\in S_2(\Gamma_0(36))$ be the new form associated to the elliptic curve $E/\BQ$.

\begin{cor}\label{cor:3}
The Heegner point $z_n$
 satisfies the following height formula
\[\frac{L'(1,E^{(n)})L(1,E^{(n^{-1})})}{\Omega^{(n)}\Omega^{(n^{-1})}}=\frac{1}{27}\widehat{h}_\BQ(z_n).\]
\end{cor}
\begin{proof}
The explicit formula in  Theorem \ref{gz-sharp} implies the following identity
\[L'(1,E,\chi_n) = \frac{(8\pi^2)(\phi,\phi)_{U_0(36)}}
    {\sqrt{3}\cdot 3N}\cdot
  \frac{\langle P_{\chi_n}(f),P_{\chi_n^{-1}}(f)\rangle_{K,K}}{(f,f)_{\CR^\times}},\]
where
\[P_{\chi_n}(f) = \sum_{t \in \widehat{K}^\times/K^\times\widehat{\CO}_{3N}^\times} f(P_0)^{\sigma_t}\chi_n(t).\]
The order $\CR$ chosen in the introduction is conjugate to $R_0(36)$ outside the place $2$. Thus
\[(f,f)_{\CR^\times} = \frac{\Vol(X_{\CR^\times})}{\Vol(X_0(36))}\deg f
= \frac{\Vol(U_0(36)_2)}{\Vol(\CR_2^\times)} = \frac{1}{3}.\]

\par Let $L = K(\sqrt[3]{n})$ and $z=\RTr_{H_{6N}/L}f(P_0)$.
When $3\nmid \sum_1^k\epsilon_i$, $\Gal(L/K)=\langle \sigma_{\omega_{3}}\rangle$. Then
\begin{eqnarray*}
\langle P_{\chi_n}(f),P_{\chi_n^{-1}}(f)\rangle_{K,K}
&=& \frac{1}{9}\left\langle\sum_{\sigma\in \Gal(L/K)} z^{\sigma}\chi_n(\sigma),
\sum_{\sigma\in \Gal(L/K)} z^{\sigma}\chi_n^{-1}(\sigma)\right\rangle_{K,K}\\
&=&\frac{1}{3}(\langle z,z\rangle_{K,K}-\langle z^{\sigma_{\omega_3}},z\rangle_{K,K} ),
\end{eqnarray*}
By Corollary $\ref{co:1}$, $z^{\sigma_{\omega_3}}\equiv [\omega]z \mod $ torsions,
\[\langle z^{\sigma_{\omega_3}},z\rangle_{K,K}=
  \frac{1}{2}\left (\widehat{h}_K([1+\omega]z)-\widehat{h}_K([\omega]z)-\widehat{h}_K(z)\right )=-\frac{1}{2}\langle z,z\rangle_{K,K}\]
 and therefore
\[\langle P_{\chi_n}(f),P_{\chi_n^{-1}}(f)\rangle_{K,K}=
  \frac{1}{2}\widehat{h}_K(z)=\widehat{h}_\BQ(z). \]
While we note that $z_n=3z$,
\[L'(1,E,\chi_n) =  \frac{8\sqrt{3}\pi^2(\phi,\phi)_{U_0(36)}}{27 N}\widehat{h}_\BQ(z_n).\]
Since
\[\Omega^{(n)}\Omega^{(n^{-1})}=\Omega^2/N=\frac{8\sqrt{3}\pi^2(\phi, \phi)_{U_0(36)}}{N},\]
we have
\[\frac{L'(1,E^{(n)})L(1,E^{(n^{-1})})}{\Omega^{(n)}\Omega^{(n^{-1})}} =  \frac{1}{27}\widehat{h}_\BQ(z_n).\]
\end{proof}

\par Now assume $k=1$. Let $p\equiv 2,5\mod 9$ be a prime number and $n=p^*$. The point $y_0=\RTr_{H_{3p}/K(\sqrt[3]{n})}(P_0-T)\in E(K(\sqrt[3]{n}))$ is of infinite order and
satisfies $y_0^{\sigma_{\omega_3}}=[\omega]y_0+t,$
with $t\in E(\BQ)[\sqrt{-3}]$ non-zero. Then $z_1=\sqrt{-3}y_0$ is a point of infinite order in $E(K(\sqrt[3]{n}))^{\chi_n}$.
\begin{lem} \label{le:a}
The point $z_1$ is not divisible by $\sqrt{-3}$ in $E(K(\sqrt[3]{n}))^{\chi_n}$.
\end{lem}
\begin{proof}
Assume $z_1=\sqrt{-3}Q+T$ with $Q\in E(K(\sqrt[3]{n}))^\chi$ and $T\in E(K(\sqrt[3]{n}))^{\chi_n}_\tor$. Suppose that $C$ is a large integer prime to $3$. Then
\[\sqrt{-3}C(y_0-Q)=CT\]
lies in $E(K(\sqrt[3]{n}))[3]$. Since $E(K(\sqrt[3]{n}))[3^\infty]=E(K)[3]$, $C(y_0-Q)$ lies in $E(K)[3]$ and assume $Cy_0=CQ+R$ with $R\in E(K)[3]$. Taking the Galois action of $\sigma_{\omega_3}$, we obtain
\[Ct=[1-\omega]R=0,\]
which is a contradiction.
\end{proof}

\begin{proof}[Proof of Theorem $\ref{th:2}$]
We showed that $\widehat{h}_\BQ(z_n)\neq 0$. Therefore $L'(1,E,\chi_n)=L'(1,E^{(n)})L(1,E^{(n^{-1})})\neq 0$ and in particular we have $\ord_{s=1}L(s,E^{(n)})=1\textrm{ and }\ord_{s=1}L(s,E^{(n^{-1})})=0$. By results of Gross-Zagier and Kolyvagin, we know $\rank_\BZ E^{(n)}(\BQ)=1$ and $\Sha(E^{(n)}/\BQ)$ is finite and that $\rank_\BZ E^{(n^{-1})}(\BQ)=0$ and $\Sha(E^{(n^{-1})}/\BQ)$ is finite. Note that Tamagawa numbers $c_v(E^{(n)})=c_v(E^{(n^{-1})})=3$ for $v\mid 2p$ and $c_v(E^{(n)})=c_v(E^{(n^{-1})})=1$ for other places $v$, and that $E^{(n)}(\BQ)\simeq \BZ/3\BZ$ and $E^{(n^{-1})}(\BQ)\simeq \BZ/3\BZ$. Then the BSD conjecture predicts
\[L'(1, E^{(n)})/\Omega^{(n)}=\# \Sha(E^{(n)}) \widehat{h}_\BQ(P),\leqno{\mathrm{BSD(n)}} \]
where $P$ is the generator of the free part of $E^{(n)}(\BQ)$, and
\[L(1, E^{(n^{-1})})/\Omega^{(n^{-1})}=\# \Sha(E^{(n^{-1})}).\leqno{\mathrm{BSD(n^{-1})}}\]
Then $\mathrm{BSD(n)\cdot BSD(n^{-1})}$ is
$$L'(1, E,\chi_n)/\Omega^{(n)}\Omega^{(n^{-1})}=\#\Sha(E^{(n^{-1})})\cdot
\#\Sha(E^{(n)})\cdot \widehat{h}_\BQ(P).$$

The descent method tells that
\[\dim_{\BF_3}\Sel_3(E^{(n)})/E^{(n)}{(\BQ)}_\tor\leq 1\textrm{ and }\dim_{\BF_3}\Sel_3(E^{(n^{-1})})/E^{(n^{-1})}{(\BQ)}_\tor= 0.\]
Hence $\Sha(E^{(n)})[3^\infty]=\Sha(E^{(n^{-1})})[3^\infty]=0$. Note $z_n=9y_0$. By the explicit Gross-Zagier formula in Corollary $\ref{cor:3}$ and the triviality of $\Sha(E^{(n)})[3^\infty]$ and $\Sha(E^{(n^{-1})})[3^\infty]$, in order to prove the $3$-part of $\mathrm{BSD(n)\cdot BSD(n^{-1})}$, it suffices to show that $\widehat{h}_\BQ(P)=u\widehat{h}_\BQ(z_1)$ with $u\in \BZ_3^\times\cap \BQ$.
\par Since the free part of $E^{(n)}(K)$ has $\CO_K$-rank $1$, the ratio $\widehat{h}_\BQ(P)/\widehat{h}_\BQ(z_1)$ is a rational number. Recall there is an isomorphism $\phi: E^{(n)}(K)\simeq E(K(\sqrt[3]{n}))^{\chi_n}$ such that $\phi(\overline {R})=\overline{\phi(R)}$ for any point $R\in E^{(n)}(\overline{\BQ})$. By Proposition $\ref{pro:5}$, $\overline{z_1}=z_1$ and hence $\phi^{-1}(z_1)$ is rational over $\BQ$. Since $z_1$ is not divisible by $\sqrt{-3}$ in $E(K(\sqrt[3]{n}))^{\chi_n}$, $\phi^{-1}(z_1)\equiv tP \mod $ torsions with some $t\in \BZ_3^\times$. Then
$\widehat{h}_\BQ(P)=u\widehat{h}_\BQ(z_1)$ with some $u\in \BZ_3^\times$.
\end{proof}


\begin{thebibliography}{10}


\bibitem{Birch} Birch, B.J., {\em Elliptic curves and modular functions}, Symposia Mathematica, Indam Rome 1968/1969, vol. 4, pp27-32. London: Academic Press (1970).

\bibitem{BSD} Birch, B.J. and Swinnerton-Dyer, H.P.,
{\em Notes on elliptic curves (II)}, J. Reine Angew. Math, 218
(1965), 79-108.

\bibitem{CST} L. Cai, J. Shu and Y. Tian {\em Explicit Gross-Zagier formular and Waldspurger formula}, 2014.

\bibitem{Cass} W. Casselman, {\em The restriction of a representation of
  $\GL_2(k)$ to $\GL_2(\CO)$}, Math. Ann. 206 (1973), 311-318.


\bibitem{CLTZ} John Coates, Yongxiong Li, Ye Tian, and Shuai Zhai,
{\em Quadratic Twists of Elliptic Curves}, preprint.


\bibitem{DV}
  S. Dasgupta and J. Voight {\em Heegner points and Sylvester's conjecture},
  Arithmetic geometry, 91-102, Clay Math. Proc., 8, Amer. Math. Soc., Providence, RI, 2009.


\bibitem{Gross1} B. Gross, {\em Heights and the special values of
L-series}. Canadian Mathematical Society, Conference Proceedings,
Volume 7 (1987).

\bibitem{Gross} B. Gross, {\em Local orders, root
numbers, and modular curves},  Americal Journal of Mathematics,
110(1988), 1153-1182.


\bibitem{GJ} S. Gelbart and H. Jacquet, {\em
		A relation between automorphic representations of
	$\GL(2)$ and $\GL(3)$},  Ann. Sci. \'Ecole Norm. Sup.
	(4) 11 (1978), no. 4, 471¨C542.

\bibitem{GP} B. Gross and D. Prasad, {\em Test vectors for linear forms},
Math. Ann, 291(1991), 343-355.

\bibitem{GZ} B. Gross and D. Zagier,  {\em Heegner points and derivatives of
$L$-series.} Invent. Math. 84 (1986), no. 2, 225--320.



\bibitem{Heegner} Heegner, K. {\em Diophantische analysis und modulfunktionen}. Math. Z. 56, 227-253 (1952).



\bibitem{Kobayashi} Shinichi Kobayshi, {\em The $p$-adic Gross-Zagier
formula for elliptic curves at supersingular primes}, Preprint.

\bibitem{K}
V.A.Kolyvagain, {\em Euler system}, The Grothendieck Festschrift.
Prog. in ath., Boston, Birkhauser (1990).

\bibitem{K1}
V.A.Kolyvagain, {\em Finiteness of $E(\BQ)$ and $\Sha(E,\BQ)$ for a
subclass of Weil curves}, Math.USSR Izvestiya,Vol. 32 (1989), No. 3.

\bibitem{Lieman} D. Lieman, {\em Nonvanishing of L-series associated to cubic twists of elliptic curves,} Ann. of Math. 140 (1994), 81???08.


\bibitem{Miyake} Toshitsune Miyake, {\em Modular Forms},
Springer Monographs in Mathematics, 1989.


\bibitem{Monsky}  P. Monsky, {\em Mock Heegner Points and Congruent Numbers}, Math. Z. 204, 45-68 (1990).

\bibitem{Nek}  Jan Nekov\'{a}\u{r}, {\em The Euler system method for
CM points on Shimura curves}, In: L-functions and Galois
representations, (Durham, July 2004), LMS Lecture Note Series 320,
Cambridge Univ. Press, 2007, pp. 471 - 547.

\bibitem{Perrin-Riou} Perrin-Riou, B., {\em  Points de Heegner et d¡äeriv¡äees de fonctions
L p-adiques}, Invent. Math. 89 (1987), no. 3, pp. 455-510.

\bibitem{Prasad} Prasad, Dipendra, {\em Some applications of seesaw
duality to braching laws}, Math.Ann. 304, 1-20 (1996).

\bibitem{Sat} H. Saito, {\em On Tunnell's formula for characters of $\GL(2)$},
  Compositio Math. 85 (1993), 99-108.
\bibitem{Satge} P. Satge, {\em Un analogue du calcul de Heegner},
Inv. Math. 87 (1987), 425?439.


\bibitem{Syl} J. J. Sylvester, {\em On certain ternary cubic-form equations}, American Journal of Math. Vol 2, No 4. 1879, pp. 357-393.

\bibitem{Tate} J. Tate, {\em Number Theoretic Background} in Automorphic Forms, Representations
  and $L$-functions, Proc. Symp. in Pure Math. XXXIII Part 2(1979), 3-26.

\bibitem{Tian1} Ye Tian, {\em Congruent Numbers with many prime
factoirs}, PNAS, Vol 109, no. 52. 21256-21258.

\bibitem{Tian2} Ye Tian, {\em Congruent Numbers and Heegner Points},
preprint.

\bibitem{Tian3} Ye Tian, {\em Congruent Number Problem},
will appear in proceeding of ICCM 2013.

\bibitem{TYZ} Ye Tian, Xinyi Yuan, and Shouwu Zhang,
{\em Genus Periods, Genus Points, and Congruent Number Problem},
preprint.

\bibitem{Tunnell1} J. Tunnell,
{\em On the local Langlands conjecture for $\GL(2)$}, Inv. Math. 46
(1978), 179-200.

\bibitem{Tunnell2} J. Tunnell, {\em Local $\epsilon$-factors and
characters of $\GL(2)$}, Amer. J. Math. 105 (1983), 1277- 1307. [13]

\bibitem{LNM800} Marie-France Vign\'{e}ras,
{\em Arithm\'{e}tique des Alg\`{e}bres de Quaternions}, Lecture Note
of Mathematics, 800. Springer, Berlin, 1980. vii+169 pp.

\bibitem{YZZ}  X.Yuan, S. Zhang, and W. Zhang,  {\em The Gross-Zagier formula on Shimura
Curves}, Annals of Mathematics Studies Number 184, 2012.




\end{thebibliography}
\end{document}